\documentclass[12pt]{article}
\input xy
\xyoption{all}
\usepackage{amsthm}
\usepackage{amsfonts,amsmath}
\usepackage{lscape}
\usepackage{amsfonts,latexsym}
\usepackage{pstricks}
\usepackage{pst-grad}
\input epsf

\def\tp{\otimes} 
\def\modn{\hspace{0.2cm} (\mbox{mod} \hspace{0.2cm} n)} 
\def\es{\vspace{0.2cm}} 

\def\d{\delta}
\def\dn{\bar{\delta}}
\newcommand\floor[1]{\left\lfloor #1 \right\rfloor}

\def\N{\mathbb{N}}
\def\Z{\mathbb{Z}}

\def\R{\mathbb{R}}
\def\C{\mathbb{C}}

\newtheorem{thm}{Theorem}[section]
\newtheorem{prop}[thm]{Proposition}

\theoremstyle{definition}

\newtheorem{lem}[thm]{Lemma}
\newtheorem{conj}[thm]{Conjecture}

\renewcommand{\proof}[1]{\noindent {\bf Proof.} 
        #1\begin{flushright}$\Box$\end{flushright}}

\theoremstyle{remark}

\newtheorem*{rmk}{Remark}

\title{Solutions of the Yang-Baxter equation: descendants of the six-vertex model from the Drinfeld doubles of dihedral group algebras}
\author{
P.E. Finch\footnote{Present address: Institut f\"{u}r Theoretische Physik, Leibniz Universit\"{a}t Hannover, Appelstra\ss e 2, 30167 Hannover, Germany}, 
K.A. Dancer, P.S. Isaac and J. Links\\
Centre for Mathematical Physics, \\
School of Mathematics and Physics, \\ 
The University of Queensland, 4072, \\
Australia.}

\date{Keywords: Quasi-triangular Hopf algebras, Yang--Baxter equation, Drinfeld double of finite groups}

\begin{document}
\maketitle

\begin{abstract}
\noindent 
The representation theory of the Drinfeld doubles of dihedral groups is used to solve the Yang-Baxter equation. Use of the $2$-dimensional representations recovers the six-vertex model solution. Solutions in arbitrary dimensions, which are viewed as descendants of the six-vertex model case,  are then obtained using tensor product graph methods which were originally formulated for quantum algebras. Connections with the Fateev-Zamolodchikov model are discussed.  
\end{abstract}

\vfil\eject
\tableofcontents

\section{Introduction}

The Yang-Baxter equation in its spectral parameter-dependent form has  well-known 
applications in a variety of areas of mathematical physics, including exactly 
solvable classical models of two-dimensional statistical mechanics and 
integrable quantum systems. 
The origins of the field trace back to the influential works by McGuire concerning  
many-body systems with delta function interactions \cite{jim}, Yang's solution for interacting fermions \cite{cny}, and Baxter's solution of the eight-vertex model \cite{rodney}. 
Since then numerous solutions have been obtained including those leading to  the Perk-Schultz model \cite{perkschultz}, the Andrews-Baxter-Forrester model \cite{abf}, the chiral Potts model \cite{cp1,cp2} and the Hubbard model \cite{shastry,usw} to name a few.

One systematic method for solving the 
Yang-Baxter equation was found with the advent of affine quantum algebras, which 
are deformations of the universal enveloping algebras of affine classical Lie 
algebras \cite{Drinfeld1986,jimbo86}.
Through the Drinfeld double construction \cite{Drinfeld1986}, the affine quantum algebras can be seen to belong to a larger class of algebras which 
are referred to as {\it quasi-triangular Hopf} algebras. The double construction provides a 
universal prescription for constructing a quasi-triangular Hopf algebra, 
denoted $D(H)$, from any Hopf algebra $H$ and its dual algebra $H^*$. Via this 
construction there exists a canonical element in the tensor product algebra of 
the Drinfeld double, $R\in D(H)\otimes D(H)$, known as the universal $R$-matrix. 
This element provides an algebraic solution of the Yang-Baxter equation. 
Thus for each matrix representation of $D(H)$, a matrix solution of the 
Yang-Baxter equation is automatically obtained. Solutions of the Yang-Baxter 
equation with spectral parameter arise whenever matrix representations of $D(H)$ 
can be parameterised by one or more continuous variables. In the context of affine 
quantum algebras, the spectral parameter arises naturally by consideration of an 
evaluation homomorphism from the affine quantum algebra to its non-affine 
counterpart. In this manner the representations of the affine quantum algebras
which give solutions with spectral parameter are the loop 
representations. In principle the evaluation homomorphism provides a route to 
algebraic solutions of the Yang-Baxter equation with spectral parameter. However 
it is a technically challenging task to undertake and has only been implemented 
explicitly for some low rank cases \cite{zg94,bgz95}.

Often, a more tractable approach is to start with a constant matrix solution of the 
Yang-Baxter equation and then ask the question of whether a spectral 
parameter-dependent generalisation exists, a procedure which is commonly referred to as 
{\it Baxterisation}. For example, in cases where the solution of the 
Yang-Baxter equation without spectral parameter affords a representation of the 
Hecke algebra, there exists a procedure for constructing the 
associated spectral parameter-dependent solution \cite{j90}. This also applies in 
cases where the constant solution gives rise to a representation of the 
Birman-Wenzl-Murakami algebra \cite{cgx91}. For solutions associated with representations 
of quantum algebras, a general prescription for Baxterising constant solutions of 
the Yang-Baxter equation (subject to satisfying certain conditions) was developed 
in \cite{zgb91,dgz94,dgz96} using the notion of the {\it tensor product graph} method. 

An alternative route to  
building solutions of the Yang-Baxter equation is in the framework of descendants.      
The most notable instances were found by Bazhanov and Stroganov 
\cite{BazhStrog1990}, who showed that the chiral Potts model is a descendant of the 
six-vertex model, and Hasegawa and Yamada \cite{HaseYamada1990}, who showed that the 
Kashiwara-Miwa model is a descendant of the zero-field eight-vertex model. In both 
cases the connection was found by constructing an appropriate algebraic structure 
through an $L$-operator. In \cite{BazhStrog1990} the algebraic structure is 
closely related to that of the quantum algebra $U_q(sl(2))$ with $q$ a root of unity, and in  
\cite{HaseYamada1990} the $L$-operator is expressed in terms of the Sklyanin algebra 
\cite{sklyanin}. In this approach use of the $L$-operator allows for the construction of 
descendants through higher dimensional representations of the algebra. It is known 
\cite{Baxter2004} that the six-vertex and the zero-field eight-vertex models intersect 
at the zero-field six-vertex model whose descendant is the Fateev-Zamolodchikov model 
\cite{FateevZam1982b}, which is precisely the intersection between the chiral Potts and
 Kashiwara-Miwa models.

Another principal class of quasi-triangular 
Hopf algebras is the class of Drinfeld doubles of finite group algebras.  These algebras have recently received attention in relation 
to the description of anyonic symmetries in quantum systems, 
specifically in terms of the braiding properties of anyonic quasiparticle excitations 
which are described by constant solutions of the Yang-Baxter equation \cite{wildb98,k03,m03,nssfd08}. However these algebras 
these have not been investigated in the context of solving the Yang-Baxter equation 
with spectral parameter, apart from some preliminary investigations reported in \cite{DIL2006,dl09}.
Arguably the simplest family of Drinfeld double algebras for finite groups is that 
associated with the dihedral groups, denoted $D_n$. The Drinfeld doubles of these 
algebras, denoted $D(D_n)$, are finite-dimensional quasi-triangular Hopf algebras with 
a finite number of irreducible representations. For $D(D_n)$ 
 all the irreducible representations are known \cite{DIL2006}. In the 
case where $n$ is odd, the irreducible representations can only have dimensions 1, 2 or 
$n$, whereas for $n$ even the dimensions are 1, 2 or $n/2$. The constant $R$-matrices 
associated with the 2-dimensional representations can be Baxterised to reproduce the 
six-vertex model at roots of unity \cite{DIL2006}. 

Our goal below is to obtain spectral 
parameter-dependent solutions of the Yang-Baxter equation associated with the
$n$-dimensional ($n$ odd) and $n/2$-dimensional ($n$ even) representations. The general 
strategy we follow is to view these solutions as descendants of the six-vertex model, 
by firstly determining the appropriate $L$-operators. Once this is achieved, the technical 
aspects for determining the explicit form of the descendants is accommodated by adapting 
the tensor product graph method developed for quantum algebras \cite{zgb91,dgz94,dgz96} 
to the present case. It turns out that solutions we obtain for $n$ odd are limiting cases of the Fateev-Zamolodchikov models \cite{BazhPerk2009}. However we stress that the anticipated $U_q(sl(2))$ symmetry with $q$ a root of unity degenerates in this limit, with the $D(D_n)$ symmetry emerging.

\section{Preliminaries}


\noindent
We first define some notation which is used extensively throughout this article.
We often work
modulo $n \in \N$, and hence define the following map: given $a \in \Z$ then $0 \leq \overline{a} \leq n-1$ is the integer 
which satisfies $a \equiv \overline{a} \modn$. We also use the following two delta functions:

\begin{equation*}
\d^j_i = 
	\begin{cases}
		1, & ~~ i = j, \\
		0, & ~~ i \neq j,
	\end{cases}
~~~~
\dn_i^j = 
	\begin{cases}
		1, & ~~ i \equiv j \modn, \\
		0, & ~~ i \not\equiv j \modn.
	\end{cases}
\end{equation*}


\noindent
We use $e_{i,j}$ to denote an elementary matrix in 
$M_{n \times n}(\C)$ whose indices are considered modulo $n$. We adopt the convention that $e_{0,0}$ corresponds to the matrix with an entry in the $n$th row and $n$th column.
These matrices obey 
the relation
$$ e_{i,j}e_{k,l} = \dn_{j}^{k}e_{i,l}. $$

\noindent
We use the following convention for the product symbol:
$$ \prod_{i=j}^{k} a_{i} = \left\{ \begin{array}{cc} 1, & k < j, \\ a_{j}a_{j+1}...a_{k}, & k \geq j. \end{array} \right. $$

\subsection{The Yang-Baxter equation with spectral parameter}
We will use three variants of the Yang-Baxter equation (YBE), which is a non-linear matrix equation in $\mbox{End }(V \tp V \tp V)$ for some vector space $V$. The first 
of these forms we call the constant YBE, given by
$$ R_{12}R_{13}R_{23} = R_{23}R_{13}R_{12}, $$
where the subscripts of $R$ refer to which vector spaces the operator is 
acting upon. That is, given
$$ R = \sum_{i} a_{i} \tp b_{i} \in \mbox{End } (V \tp V), $$
we have
$$ R_{12} = \sum_{i} a_{i} \tp b_{i} \tp I, \hspace{1cm} R_{13} = \sum_{i} a_{i} \tp I \tp b_{i}, \hspace{1cm} \mbox{etc.} $$
Here $I$ denotes the identity matrix. The second form is known as the spectral parameter-dependent YBE and is given by

\begin{equation} \label{YBE1}
 R_{12}(x)R_{13}(xy)R_{23}(y) = R_{23}(y)R_{13}(xy)R_{12}(x),
\end{equation}

\noindent
where $x,y \in \C$. 
Given a parameter-dependent solution, $R(z)$, constant solutions are 
recovered when we take the limits $z \rightarrow 0,1$ and $\infty$. We 
refer to invertible solutions of either form as {\it $R$-matrices}. 
As this article deals primarily with 
the parameter-dependent YBE we shall henceforth refer to it simply as the 
Yang-Baxter equation.

\

\noindent
Now suppose $r(z)\in \mbox{End }(V \tp V)$ and $R(z)\in \mbox{End }(W \tp W)$ are $R$-matrices  
with $\dim V < \dim W$. We describe $R(z)$ as a {\it descendant} of $r(z)$ provided 
that there exists a non-trivial invertible operator $L(z) \in \mbox{End }(V \tp W)$ (referred to as an {\it $L$-operator}) 
which satisfies
\begin{equation}
	\label{eqnrLL}
	r_{12}(xy^{-1})L_{13}(x)L_{23}(y) = L_{23}(y)L_{13}(x)r_{12}(xy^{-1}), \quad \forall x,y \in \C,
\end{equation}
and
\begin{equation}
  \label{eqnLLR} 
  L_{12}(x)L_{13}(y)R_{23}(x^{-1}y) = R_{23}(x^{-1}y)L_{13}(y)L_{12}(x), \quad \forall x,y \in \mathbb{C}.
 \end{equation}

In this article we start with an $R$-matrix associated with the zero-field 
six-vertex model and use the Drinfeld doubles of dihedral groups to generate 
$L$-operators and descendants. 

The third variant of the YBE used in this paper is
$$ \check{R}_{12}(x)\check{R}_{23}(xy)\check{R}_{12}(y) = \check{R}_{23}(y)\check{R}_{12}(xy)\check{R}_{23}(x). $$
Given a matrix solution $R(z)$ to Equation \eqref{YBE1}, a solution to this form of the YBE is given by
$$
\check{R}(z) = P\, R(z),
$$
where 
$P$ is the permutation operator, i.e.
$$ P = \sum_{i,j} e_{i,j} \tp e_{j,i}. $$

\subsection{Drinfeld doubles of finite groups}

Next we recall relevant results from the Drinfeld double construction applied to finite group algebras, following \cite{DPR1990,Gould1993}. Given a group $G$, with identity $e$, we consider the algebra $\C G$, whose basis 
vectors are the elements of the group. The multiplication and unit of $\C G$ 
are inherited from the group in the natural way. We equip $\C G$ with a coproduct, counit 
and antipode defined respectively by
$$ \Delta(g) = g \tp g, \hspace{1cm} \epsilon(g) = 1 \hspace{0.7cm} \mbox{and} \hspace{0.7cm} \gamma(g) = g^{-1}, \hspace{1cm} \forall g \in G. $$
With these maps $\C G$ becomes a Hopf algebra. We next consider the dual 
space of $\C G$, 
$$ (\C G)_{o} = \C \{g^{*} | g \in G \}. $$
The multiplication and unit are given, respectively, by
$$ g^{*}h^{*} = \d_{g}^{h} \hspace{1cm} \mbox{and} \hspace{1cm} u(1) = \sum_{g\in G}g^{*}. $$
The costructure and antipode are defined by
$$ \Delta(g^{*}) = \sum_{h\in G} (hg)^{*} \tp (h^{-1})^{*}, \hspace{1cm} \epsilon(g^{*}) = \d_{g}^{e} \hspace{0.7cm} \mbox{and} \hspace{0.7cm} \gamma(g^{*}) = (g^{-1})^{*}, \hspace{1cm} \forall g \in G. $$
Under these maps $(\C G)_{o}$ is also a Hopf algebra. Using the dual and the 
original algebra we can construct the Drinfeld double, 
$$ D(G) = \C \{gh^{*} | g,h \in G \}. $$
We impose the relation
$$ h^{*}g = g(g^{-1}hg)^{*} $$
and adopt the required remaining structure from $\C G$ and $(\C G)_{o}$. 
Furthermore it is known that $D(G)$ is a quasi-triangular Hopf algebra, containing the canonical element
$$ \mathcal{R} = \sum_{g \in G} g \tp g^{*}. $$
This element satisfies the following relations:
$$ \begin{array}{rcll}
  \mathcal{R} \Delta(a) & = & \Delta^{T}(a) \mathcal{R}, & \forall a \in D(G), \\
  (\Delta \tp \mbox{id}) \mathcal{R} & = &  \mathcal{R}_{13}\mathcal{R}_{23}, \\
  (\mbox{id} \tp \Delta) \mathcal{R} & = &  \mathcal{R}_{13}\mathcal{R}_{12}, \\
\end{array} $$
where $\Delta^T$ denotes the opposite coproduct.  Consequently matrix representations of $\mathcal{R}$ provide solutions of the constant Yang-Baxter equation. 

In this paper we use the dihedral group $D_n$, which is the symmetry group of a regular polygon with $n$ vertices.  That is,
$$ D_{n} = \{\sigma, \tau | \sigma^{n} = \tau^{2} = e, \sigma \tau \sigma = \tau \}, $$
and has order $2n$. For ease of calculation we divide the description of the representation theory into different cases. 

\subsection{Representations of $D(D_{n})$} \label{r}

\subsubsection{The case when $n$ is odd} 
We first consider $D(D_{n})$ for the case where $n$ is odd. As stated in \cite{DIL2006} 
the representations of the double of a group are naturally partitioned by the 
conjugacy classes of the group. For these representations we consider $w \in \C$ 
to be a primitive $n$th root of unity. 
Then the irreducible representations (irreps) are given in Table \ref{oddreps} below:

\begin{table}[ht]

$$
\begin{array}{|c|c|c|c|c|} \hline
 \mbox{Irrep }\pi & \mbox{Constraints} & \pi(\sigma) & \pi(\tau) & \pi(g^*),\; g \in D_n \\ \hline
	\pi_1^\pm & & 1 & \pm 1 & \delta^e_g \\[2pt]
	\pi_2^{(0,k)} & 1 \leq k \leq \frac{n-1}{2} & \begin{pmatrix} w^k & 0 \\ 0 & w^{-k} \end{pmatrix}
		& \begin{pmatrix} 0 & 1 \\ 1 & 0 \end{pmatrix} & \begin{pmatrix} \delta^e_g  & 0 \\ 0 & \delta^e_g  \end{pmatrix} \\
		[12pt]
	\pi_2^{(l,k)} & \begin{array}{c} 1 \leq l \leq \frac{n-1}{2}, \\ 0 \leq k \leq n-1 \end{array} & 
		\begin{pmatrix} w^k & 0 \\ 0 & w^{-k} \end{pmatrix} & \begin{pmatrix} 0 & 1 \\ 1 & 0 \end{pmatrix} & 
		\begin{pmatrix} \delta^{\sigma^l}_g & 0 \\ 0 & \delta^{\sigma^{-l}}_g \end{pmatrix} \\ [12pt]
	\pi_n^\pm & & \sum_{i=1}^n e_{i+1, i} & \pm \sum_{i=1}^n e_{i, 2-i} & \delta^{\sigma^{2j} \tau}_g e_{j+1, j+1} \\
[2pt] \hline
\end{array}$$ \vspace{-6mm}
\caption{\label{oddreps} The irreps of $D(D_n)$ when $n$ is odd.} 
\end{table}

\noindent Here $\pi_1^\pm, \pi_2^{(a,b)}$ and $\pi_n^\pm$ have dimensions 1, 2 and $n$ respectively, 
and their associated modules are denoted $V_1^\pm, V_2^{(a,b)}$ and $V_n^\pm$.

Also associated with these representations are the solutions of the Yang-Baxter 
equation arising from the canonical element. If we apply a two-dimensional irrep 
we find the $R$-matrix
$$ (\pi_{2}^{(l,k)} \tp \pi_{2}^{(l,k)}) \mathcal{R} = 
 \left(
	\begin{array}{cccc}
		w^{kl} & 0 & 0 & 0\\
		0 & w^{-kl} & 0 & 0\\
		0 & 0 & w^{-kl} & 0\\
		0 & 0 & 0 & w^{kl}
	\end{array}
\right). $$
It was shown in \cite{DIL2006} that this constant solution leads to the 
parameter-dependent solution
\begin{equation}
r(z) = 
	\left(
	\begin{array}{cccc}
		w^{2kl}z^{-1} - w^{-2kl}z & 0 & 0 & 0\\
		0 & z^{-1} - z & w^{2kl} - w^{-2kl} & 0\\
		0 & w^{2kl} - w^{-2kl} & z^{-1} - z & 0\\
		0 & 0 & 0 & w^{2kl}z^{-1} - w^{-2kl}z
	\end{array}
\right). \label{rznodd}
\end{equation}
This $R$-matrix corresponds to the six-vertex model with zero-field. Also of interest 
is the representation of the canonical element using the $n$-dimensional irreps,
$$ (\pi_{n}^{\pm} \tp \pi_{n}^{\pm}) \mathcal{R} = \pm \sum_{i,j=0}^{n-1}  e_{i+j,i-j} \tp e_{i,i}. $$

In our construction of the descendant of $r(z)$ we will use a linear combination 
of projection operators. It is known \cite{Gould1993} that there exist operators 
which project $D(D_{n})$ onto its ideals. These projection operators are 
defined by
$$E_{\alpha} = \frac{d[\alpha]}{|G|} \sum_{g,h \in G} 
\chi_{\alpha}(h^{*}g^{-1})gh^{*}, $$
where $\alpha$ is an irrep, $d[\alpha]$ its dimension and $\chi_{\alpha}$ is the 
group character defined by
$$
\chi_\alpha(a) = \mbox{tr } \pi_\alpha(a),\ \ \forall a\in D(D_n). 
$$
We also consider the projection operators
\begin{equation}
	p_n^{\alpha} = (\pi_{n}^\pm \tp \pi_{n}^\pm) \Delta (E_{\alpha}). 
\label{eqnalgebraicprojection2}
\end{equation}
These operators project from $\pi_{n}^\pm \tp
\pi_{n}^\pm$ onto copies of the irrep associated with $\alpha$ in the
decomposition of $\pi_{n}^\pm \tp \pi_{n}^\pm$. 
Note that we consider $\pi_n^+\tp\pi_n^+$ and $\pi_n^-\tp\pi_n^-$ together
since  
$$ (\pi_{n}^{+} \tp \pi_{n}^{+})\Delta(a) = (\pi_{n}^{-} \tp \pi_{n}^{-})\Delta(a), \hspace{1cm} \forall a \in D(D_{n}). $$
We calculate $p_n^\alpha$ explicitly, using the expression given in Equation (\ref{eqnalgebraicprojection2}), obtaining
\begin{equation*}
p_n^{\alpha} 
 =  \frac{d[\alpha]}{2n} \sum_{g \in D_{n}} 
\sum_{i,j = 0}^{n-1}\chi_{\alpha}((\sigma^{2j})^{*}g^{-1}) 
\pi_{n}^\pm(g) e_{i-j,i-j} \tp \pi_{n}^\pm(g) e_{i,i}.
\end{equation*}
From this we can see that the only irreps with non-zero projection 
operators will be associated with the conjugacy classes $\{e\}$ and 
$\{\sigma^{i}, \sigma^{-i} \}$ for $ 1 \leq i \leq \frac{n-1}{2}$. 
This implies that $\pi_{n}^\pm \tp \pi_{n}^\pm$ decomposes only into one 
and two-dimensional irreps.

For convenience, we slightly modify our notation for the irreps. Instead of
using $\pi_1^+$ and $\pi_2^{(l,k)}$, we use ordered pairs corresponding only to
irreps that appear in the direct sum decomposition of $\pi_n^\pm\tp\pi_n^\pm$.
The correspondence is summarised in Table \ref{alphaodd}:

\begin{table}[ht]
$$ \begin{array}{|c|c|c|c|} \hline
  \alpha & \mbox{Irrep} & \mbox{Constraint on } a & \mbox{Constraint on } b \\ \hline
  (0,0) & \pi_{1}^{+} & & \\
  (0,b) & \pi_{2}^{(0,2b)} & & 1 \leq b \leq \floor{\frac{n-1}{4}} \\
  (0,b) & \pi_{2}^{(0,n-2b)} & & \floor{\frac{n+3}{4}} \leq b \leq \frac{n-1}{2} \\
  (a,b) & \pi_{2}^{(2a,\overline{2b})} & 1 \leq a \leq \floor{\frac{n-1}{4}} & 0 \leq b \leq n-1 \\
  (a,b) & \pi_{2}^{(n-2a,\overline{n-2b})} & \floor{\frac{n+3}{4}} \leq a \leq \frac{n-1}{2} & 0 \leq b \leq n-1 \\ \hline \end{array} $$ \vspace{-6mm}
\caption{\label{alphaodd} The ordered pairs labelling the irreps for $D(D_n)$, $n$ odd.} 
\end{table}
\noindent Here $\floor{a}$ denotes the floor of $a$.
Then the projection operator for irrep $\alpha = (a,b)$ is given by 
$$ p_n^{\alpha} = \frac{c^\alpha}{n} \sum_{i,j = 0}^{n-1}[ w^{2bj} e_{i+a+j,i+a} \tp e_{i+j,i} + w^{-2bj} e_{i-a+j,i-a} \tp e_{i+j,i}]$$
 
\noindent where

$$c^\alpha = \begin{cases}
\frac{1}{2}, &\quad \alpha = (0,0), \\
1, & \quad \alpha \neq (0,0).
\end{cases}$$

We have calculated non-zero projection operators for $(n^2-1)/2$ two-dimensional irreps and one 1-dimensional irrep.  
Using a counting argument, it is clear that these provide a complete decomposition of $\pi_n^\pm \tp \pi_n^\pm$ into irreps for $n$ odd.



\subsubsection{The case when $n$ is even \label{sub2-3-2}}
We similarly catalogue the irreps of $D(D_{n})$ for the case when $n$ is even.
We set $n=2m$ and let $w \in \C$ be a primitive $2m$th 
root of unity. 

\begin{table}[ht]
$$
\begin{array}{|c|c|c|c|c|} \hline
 \mbox{Irrep }\pi & \mbox{Constraints} & \pi(\sigma) & \pi(\tau) & \pi(g^*),\; g \in D_{2m} \\ \hline
	\pi_{1,e}^{(a,b)} & a,b \in \{0,1\} & (-1)^b & (-1)^a & \delta^e_g \\[5pt]
	\pi_{1,\sigma^m}^{(a,b)} & a, b \in \{0,1\} & (-1)^b & (-1)^a & \delta^{\sigma^m}_g \\[5pt]
	\pi_2^{(0,k)} & 1 \leq k < m & \begin{pmatrix} w^k & 0 \\ 0 & w^{-k} \end{pmatrix}
		& \begin{pmatrix} 0 & 1 \\ 1 & 0 \end{pmatrix} & \begin{pmatrix} \delta^e_g  & 0 \\ 0 & \delta^e_g  \end{pmatrix} \\
		[12pt]
	\pi_2^{(m,k)} & 1 \leq k < m & \begin{pmatrix} w^k & 0 \\ 0 & w^{-k} \end{pmatrix}
		& \begin{pmatrix} 0 & 1 \\ 1 & 0 \end{pmatrix} & 
		\begin{pmatrix} \delta^{\sigma^m}_g  & 0 \\ 0 & \delta^{\sigma^m}_g  \end{pmatrix} \\[12pt]
	\pi_2^{(l,k)} & \begin{array}{c} 1 \leq l \leq m-1, \\ 0 \leq k \leq 2m-1 \end{array} & 
		\begin{pmatrix} w^k & 0 \\ 0 & w^{-k} \end{pmatrix} & \begin{pmatrix} 0 & 1 \\ 1 & 0 \end{pmatrix} & 
		\begin{pmatrix} \delta^{\sigma^l}_g & 0 \\ 0 & \delta^{\sigma^{-l}}_g \end{pmatrix} \\ [12pt]
	\pi_{m,\tau}^{(a,b)} & a, b \in \{0,1\} & \sum_{i=1}^m (-1)^{a\delta_{i}^{1}}e_{i,i-1} & \sum_{i=1}^m 
		(-1)^{a\delta_{i}^{1}+b}e_{i,2-i} &	\delta^{\sigma^{2k}\tau}_g e_{k+1, k+1} \\ [5pt]
	\pi_{m, \sigma \tau}^{(a,b)} & a, b \in \{0,1\} & \sum_{i=1}^m (-1)^{a\delta_{i}^{1}}e_{i,i-1} & (-1)^b \sum_{i=1}^m e_{i,1-i} &
		\delta^{\sigma^{2k+1}\tau}_g e_{k+1, k+1}  \\
[2pt] \hline
\end{array}$$ \vspace{-6mm}
\caption{The irreps of $D(D_{2m})$.}
\end{table}

We again concern ourselves with the two-dimensional representations applied to 
the canonical element. We find the constant solution 
$$ (\pi_{2}^{(l,k)} \tp \pi_{2}^{(l,k)}) \mathcal{R} = 
 \left(
	\begin{array}{cccc}
		w^{kl} & 0 & 0 & 0\\
		0 & w^{-kl} & 0 & 0\\
		0 & 0 & w^{-kl} & 0\\
		0 & 0 & 0 & w^{kl},
	\end{array}
\right), $$
which leads to the parameter-dependent solution 
$$ r(z) = 
	\left(
	\begin{array}{cccc}
		w^{2kl}z^{-1} - w^{-2kl}z & 0 & 0 & 0\\
		0 & z^{-1} - z & w^{2kl} - w^{-2kl} & 0\\
		0 & w^{2kl} - w^{-2kl} & z^{-1} - z & 0\\
		0 & 0 & 0 & w^{2kl}z^{-1} - w^{-2kl}z
	\end{array}
\right). $$
This result only differs to that of $D(D_{n})$ where $n$ is odd by the possible 
choices of the root of unity. 
Note also the similarity in the representation of the canonical element, given by 
$$ (\pi_{m} \tp \pi_{m})\mathcal{R} = (-1)^{b} \sum_{i,j=0}^{m-1} e_{i+j,i-j} \tp e_{i,i}, $$
when $\pi_{m}$ is either $\pi_{m,\tau}^{(0,b)}$ or $\pi_{m,\sigma\tau}^{(0,b)}$ for $b\in \{0,1\}$.

We again use Equation $(\ref{eqnalgebraicprojection2})$ to derive our projection
operators. We first look at the case when $m=n/2$ is odd. We find that
$$ 
p_m^{\alpha} = \frac{d[\alpha]}{4m} \sum_{g \in G} \sum_{i,j=0}^{m-1} 
\chi_{\alpha}((\sigma^{2j})^{*}g^{-1}) \pi_{m}(g) e_{i-j,i-j} \tp \pi_{m}(g) 
e_{i,i}, 
$$
where $\pi_{m}$ is any one of the $m$-dimensional irreps. 
As in the case of $D(D_n)$ for odd $n$, we introduce an ordered pair notation for the irreps
that appear in the direct sum decomposition of $\pi_{m,\tau}^{(0,b)} \tp \pi_{m,\tau}^{(0,b)}$ and $\pi_{m,\sigma\tau}^{(0,b)} \tp \pi_{m,\sigma\tau}^{(0,b)}$ for $b\in \{0,1\}$.
The correspondence is summarised in Table \ref{alpha2odd}:
\begin{table}[ht]
$$ \begin{array}{|c|c|c|c|} \hline
  \alpha & \mbox{Irrep} & \mbox{Constraint on } a & \mbox{Constraint on } b \\ \hline
  (0,0) & \pi_{1,e}^{+} & & \\
  (0,b) & \pi_{2}^{(0,2b)} & & 1 \leq b \leq \frac{m-1}{2} \\
  (a,b) & \pi_{2}^{(2a,2b)} & 1 \leq a \leq \frac{m-1}{2} & 0 \leq b \leq m-1 \\ \hline 
\end{array} $$ \vspace{-6mm}
\caption{\label{alpha2odd}The ordered pairs labelling the irreps for $D(D_{2m})$, $m$ odd.}
\end{table}

%
%
%
Then the projection operator associated with irrep $\alpha = (a,b)$ is given by
$$ p_m^{\alpha} = \frac{c^\alpha}{m} \sum_{i,j=0}^{m-1}
[w^{2bj} e_{i+a+j,i+a} \tp e_{i+j,i} +  w^{-2bj} e_{i-a+j,i-a} \tp e_{i+j,i}]$$

\noindent where $$c^\alpha = \begin{cases}
\frac{1}{2}, &\quad \alpha = (0,0), \\
1, & \quad \alpha \neq (0,0).
\end{cases}$$

We note that these projection operators here match those arising in the case of $D(D_{n})$ where is $n$ odd. 

Similarly, when $m=n/2$ is even we have the projection operators
$$ p_m^{\alpha} = \frac{d[\alpha]}{4m} \sum_{g \in G} \sum_{i,j=0}^{m-1} 
\chi_{\alpha}((\sigma^{2j})^{*}g^{-1}) \pi_{m}(g) e_{i-j,i-j} \tp \pi_{m}(g) 
e_{i,i}, $$
for any $m$-dimensional irrep $\pi_m$. 
The ordered pair notation for the irreps occuring in the decomposition of 
$\pi_{m,\tau}^{(0,b)} \tp \pi_{m,\tau}^{(0,b)}$ for $b\in \{0,1\}$ is given in Table \ref{case i}, whereas 
those occuring in the decomposition of $\pi_{m,\sigma\tau}^{(0,b)} \tp \pi_{m,\sigma\tau}^{(0,b)}$ for $b\in \{0,1\}$ are given in 
Table \ref{case ii}.
\begin{table}[ht]
$$ \begin{array}{|c|c|c|c|} \hline
  \alpha & \mbox{Irrep} & \mbox{Constraint on } a & \mbox{Constraint on } b \\ \hline
  (0,b\frac{m}{2}) & \pi_{1,e}^{(0,b)} & & b \in \{0,1\} \\
  (\frac{m}{2},b\frac{m}{2}) & \pi_{1,\sigma^{m}}^{(0,b)} & & b \in \{0,1\} \\
  (0,b) & \pi_{2}^{(0,2b)} & & 1 \leq b \leq \frac{m}{2} - 1 \\
  (\frac{m}{2},b) & \pi_{2}^{(m,2b)} & & 1 \leq b \leq \frac{m}{2} - 1 \\
  (a,b) & \pi_{2}^{(2a,2b)} & 1 \leq a \leq \frac{m}{2} - 1 & 0 \leq b \leq m-1 \\ \hline \end{array} $$
\vspace{-6mm}
\caption{The ordered pairs labelling the irreps for $D(D_{2m})$, $m$ even, case $(i)$. \label{case i}}
\end{table}


\begin{table}[ht]
$$ \begin{array}{|c|c|c|c|} \hline
  \alpha & \mbox{Irrep} & \mbox{Constraint on } a & \mbox{Constraint on } b \\ \hline
  (0,b\frac{m}{2}) & \pi_{1,e}^{(b,b)} & & b \in \{0,1\} \\
  (\frac{m}{2},b\frac{m}{2}) & \pi_{1,\sigma^{m}}^{(b,b)} & & b \in \{0,1\} \\
  (0,b) & \pi_{2}^{(0,2b)} & & 1 \leq b \leq \frac{m}{2} - 1 \\
  (\frac{m}{2},b) & \pi_{2}^{(m,2b)} & & 1 \leq b \leq \frac{m}{2} - 1 \\
  (a,b) & \pi_{2}^{(2a,2b)} & 1 \leq a \leq \frac{m}{2} - 1 & 0 \leq b \leq m-1 \\ \hline \end{array} $$
\vspace{-6mm}
\caption{The ordered pairs labelling the irreps for $D(D_{2m})$, $m$ even, case $(ii)$. \label{case ii}}
\end{table}

We obtain the same projection operators irrespective of which of the two cases we consider.  In particular, the 
projection operator corresponding to the irrep $\alpha = (a,b)$ is

$$ p_m^{\alpha} = \frac{c^\alpha}{m} \sum_{i,j=0}^{m-1} 
[w^{2bj} e_{i+a+j,i+a} \tp e_{i+j,i} + w^{-2bj} e_{i-a+j,i-a} \tp  e_{i+j,i}] $$

\noindent where $$c^\alpha = \begin{cases}
\frac{1}{2}, &\quad a \in \{0,\frac{m}{2}\} \mbox{ and } b \in \{0, \frac{m}{2}\}, \\
1, & \quad \mbox{otherwise}.
\end{cases}$$

Again, a counting argument verifies that these are the only non-zero projection operators. 
As the projection operators are the same for the different $m$-dimensional irreps, we can 
without loss of generality consider only
$\pi_{m,\tau}^{(0,0)}.$

\section{Descendants associated with $D(D_{n})$ - $n$ odd}
In this section we construct a family of solutions of the Yang-Baxter equation 
using $D(D_{n})$, where $n$ is odd and $n>2$. 

\subsection{Construction of the $L$-operator\label{ssecLopDnodd}}
We begin by constructing an operator $L(z) \in {\rm End}(V_{2}^{(l,k)} \tp V_{n}^\pm)$. We observe that 
$V_{2}^{(l,k)} \tp V_n^\pm \cong V_n^+ \oplus V_n^-$, suggesting that the $L$-operator will be a linear 
combination of two intertwining operators with one parameter-dependent coefficient (after scaling). Hence we adopt the following 
ansatz for the $L$-operator: 

$$ L(z) = (\pi_{2}^{(l,k)} \tp \pi_{n}^{\pm}) [\mathcal{R} + h(z) (\mathcal{R}^{T})^{-1}] $$

\noindent where $(\mathcal{R}^{T})^{-1} = \underset{g \in G}{\sum} g^* \tp g^{-1}$. Applying the representations leads to the following expression for $L(z)$:
\begin{equation*} 
L(z) = \sum_{i=0}^{n-1}\left \{ (w^{2(i-1)k}e_{1,2} 
+ w^{-2(i-1)k}e_{2,1}) \tp e_{i,i} + h(z) \left[ e_{1,1} \tp e_{i-l,i}
 + e_{2,2} \tp e_{i+l,i} \right] \right\}.  
\end{equation*}
Performing a basis transformation which leaves $r(z)$ invariant, we obtain
$$ L(z) = \sum_{i=0}^{n-1}\left \{ (w^{2ik}e_{1,2} + w^{-2ik}e_{2,1}) \tp e_{i,i} 
+ h(z) \left[ e_{1,1} \tp e_{i-l,i} + e_{2,2} \tp e_{i+l,i} \right] \right\}. $$

Substituting $r(z)$ and $L(z)$ into Equation $(\ref{eqnrLL})$, 
we find only one constraint on $h(z)$, namely
$$ h(x)y - h(y)x = 0, $$
which has the solution
$$ h(z) = Cz, $$
for any $C \in \C$. We are free to rescale the parameter $z$ without affecting the descendants given in the next section, so without 
loss of generality we can choose $C = 1$. We have therefore shown the following:
\begin{prop}
The $L$-operator given explicitly by
\begin{eqnarray}
 L(z) = \sum_{i=0}^{n-1}\left \{ (w^{2ik}e_{1,2} + w^{-2ik}e_{2,1}) \tp e_{i,i} 
+ z \left[ e_{1,1} \tp e_{i-l,i} + e_{2,2} \tp e_{i+l,i} \right] \right\}, 
\label{lop}
\end{eqnarray}
and the $r(z)$ given in Equation (\ref{rznodd}) together satisfy Equation (\ref{eqnrLL}). 
\end{prop}
It is important to comment that, up to a basis transformation, the $L$-operator (\ref{lop}) is a particular limit 
of the general $L$-operator discussed in \cite{BazhStrog1990} in relation to the chiral Potts model. However, the 
associated $U_q(sl(2))$ structure with $q^3=1$ described in \cite{BazhStrog1990} is lost in this limit. 
Specifically,  Equation (2.9) in \cite{BazhStrog1990} does not hold. Amongst the defining relations of the generalised 
$U_q(sl(2))$ algebra this means, in particular, that  
$$\left[e,\,f\right]=0 $$ 
where $e$ and $f$ are the raising and lowering generators respectively. The above relation is indicative of the fact 
that the $U_q(sl(2))$ structure degenerates in this limit. Instead we find that the symmetry algebra of  $D(D_n)$ 
emerges when the $L$-operator is given by Equation (\ref{lop}).

\subsection{Construction of the descendants\label{props}}
\begin{figure}[ht]
\begin{center}
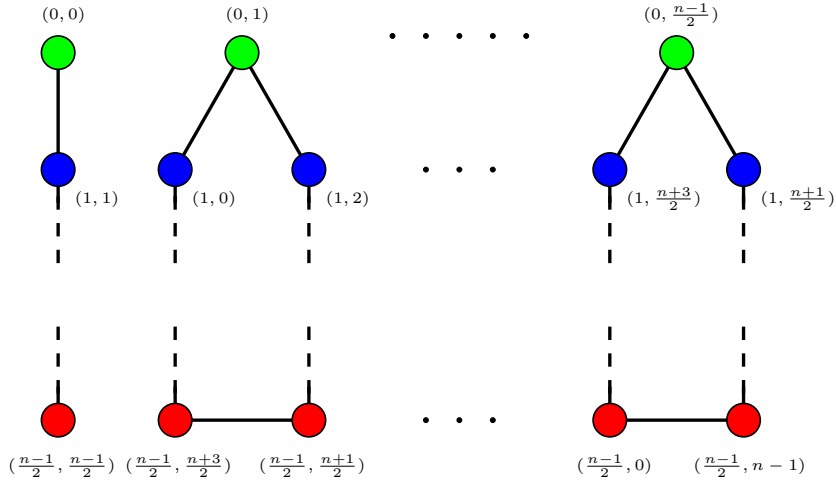
\end{center}
\caption{Tensor product graph for odd $n$ when $l = k = 1$. \label{FigTensorOdd}} 
\end{figure}

Here we put forth a predicted form of a descendant, $R(z)$, and determine the 
constraints on this form. We then impose additional constraints enforcing that 
the descendant inherits certain properties from representations of the canonical element 
from $D(D_{n})$. For convenience we shall use an alternate form of Equation $(\ref{eqnLLR})$:
\begin{equation}
	\check{R}_{23}(xy^{-1})L_{13}(x)L_{12}(y) 
= L_{13}(y)L_{12}(x)\check{R}_{23}(xy^{-1}), \label{eqnLLRc}
\end{equation}
where
$$ R(z) = P\check{R}(z)$$
and $L(z)$ is given by Equation \eqref{lop}.  As the descendants $\check{R}(z)$ commute with the action of the coproduct, 
we assume they are of the form
\begin{equation} 
\check{R}(z) = \sum_{\alpha \in S} f_{\alpha}(z) p^{\alpha}, \label{Rcheck}
\end{equation}
where $p^{\alpha}$ are the projection operators previously calculated, 
$f_{\alpha}(z)$ are continuous functions and $S$ is the set of ordered pairs 
which correspond to non-zero projection operators. For convenience we rescale our projection operators, henceforth using
\begin{equation*} 
\tilde{p}^{\alpha} 
 = \sum_{i,j = 0}^{n-1}[ w^{2jb} e_{i+a+j,i+a} \tp e_{i+j,i} + w^{-2bj} e_{i-a+j,i-a} \tp e_{i+j,i}], \quad \alpha = (a,b) \in S.
\end{equation*}
Applying 
$\tilde{p}^{\alpha}_{23}$ and $\tilde{p}^{\beta}_{23}$ to the left- and right-hand sides
respectively of Equation $(\ref{eqnLLRc})$, we find
\begin{equation} 
f_{\alpha}(xy^{-1}) \tilde{p}^{\alpha}_{23} L_{13}(x)L_{12}(y) \tilde{p}_{23}^{\beta} 
= f_{\beta}(xy^{-1}) \tilde{p}^{\alpha}_{23} L_{13}(y)L_{12}(x) \tilde{p}_{23}^{\beta}. 
\label{eqnpLLp}
\end{equation}
We calculate that
$$ 
L_{13}(x)L_{12}(y) = e_{11} \tp A(x,y) + e_{12} \tp B(x,y) + e_{21} \tp C(x,y) + e_{22} \tp D(x,y), 
$$
where
\begin{eqnarray}
	A(x,y) & = & \sum_{i,j=0}^{n-1} w^{2(j-i)k}  e_{i,i} \tp e_{j,j} + xy\ e_{i-l,i} \tp e_{j-l,j}, \nonumber\\
	B(x,y) & = & \sum_{i,j=0}^{n-1} x w^{2ik} e_{i,i} \tp e_{j-l,j} + y w^{2jk} e_{i+l,i} \tp e_{j,j}, \nonumber\\
	C(x,y) & = & \sum_{i,j=0}^{n-1} y w^{-2jk} e_{i-l,i} \tp e_{j,j} + x w^{-2ik} e_{i,i} \tp e_{j+l,j}, \nonumber\\
	D(x,y) & = & \sum_{i,j=0}^{n-1} w^{2(i-j)k} e_{i,i} \tp e_{j,j} + xy\ e_{i+l,i} \tp e_{j+l,j}. \label{ABCD}
\end{eqnarray}
Substituting these functions into Equation $(\ref{eqnpLLp})$, we find the following constraints on $f_{\alpha}(xy^{-1})$:
\begin{eqnarray}
f_{\alpha}(xy^{-1}) \tilde{p}^{\alpha} A(x,y) \tilde{p}^{\beta}	
&=&	f_{\beta}(xy^{-1}) \tilde{p}^{\alpha} A(x,y) \tilde{p}^{\beta},\label{Aconstraint} \\
f_{\alpha}(xy^{-1}) \tilde{p}^{\alpha} B(x,y) \tilde{p}^{\beta} 
&=&	f_{\beta}(xy^{-1}) \tilde{p}^{\alpha} B(y,x) \tilde{p}^{\beta},\label{Bconstraint} \\
f_{\alpha}(xy^{-1}) \tilde{p}^{\alpha} C(x,y) \tilde{p}^{\beta}	
&=&	f_{\beta}(xy^{-1}) \tilde{p}^{\alpha} C(y,x) \tilde{p}^{\beta},\label{Cconstraint} \\
f_{\alpha}(xy^{-1}) \tilde{p}^{\alpha} D(x,y) \tilde{p}^{\beta}	
&=&	f_{\beta}(xy^{-1}) \tilde{p}^{\alpha} D(x,y) \tilde{p}^{\beta}. \label{Dconstraint}
\end{eqnarray}
Note that Equations (\ref{Aconstraint},\ref{Dconstraint}) are automatically satisfied. 
The above relations provide a tensor product graph \cite{zgb91,dgz94,dgz96}. Assigning each irrep to a vertex of a graph, we connect the  
vertices labelled $\alpha$ and $\beta$ by an edge if 
$$\tilde{p}^{\alpha} \chi(x,y) \tilde{p}^{\beta}\neq 0$$	
for either $\chi = B$ or $\chi = C$. The tensor product graph for the case $l=k=1$ is depicted in Figure \ref{FigTensorOdd}. 
An edge connecting vertices $\alpha$ and $\beta$ signifies that the functions $f_\alpha(x)$ and $f_\beta(x)$ are constrained by 
Equations (\ref{Bconstraint},\ref{Cconstraint}).   

We now give a series of propositions and lemmas which result in a solution to Equation \eqref{eqnLLRc}.
\begin{prop} \label{prop3-4}
For $\check{R}(z)$ defined by Equation (\ref{Rcheck}) to be a solution of 
Equation (\ref{eqnLLRc}), the following four constraints must be satisfied:
\begin{eqnarray*}
	0 & = & \dn_{k+b}^{d} \dn_{a+l}^{c} \left[f_{(a,b)}(z)(zw^{2((-l-a)k-bl)} + 1) - f_{(c,d)}(z)(w^{2((-l-a)k-bl)} + z) \right], \\	
	0 & = & \dn_{k-b}^{d} \dn_{c}^{l-a} \left[f_{(a,b)}(z)(zw^{2((-l+a)k+bl)} + 1) - f_{(c,d)}(z)(w^{2((-l+a)k+bl)} + z) \right], \\	
	0 & = & \dn_{k+b}^{-d} \dn_{-c}^{l+a} \left[f_{(a,b)}(z)(zw^{2((-l-a)k-bl)} + 1) - f_{(c,d)}(z)(w^{2((-l-a)k-bl)} + z) \right], \\	
	0 & = & \dn_{k-b}^{-d} \dn_{c}^{a-l} \left[f_{(a,b)}(z)(zw^{2((-l+a)k+bl)} + 1) - f_{(c,d)}(z)(w^{2((-l+a)k+bl)} + z) \right], 
\end{eqnarray*}
where the pairs $(a,b)$ and $(c,d)$ belong to $S$. 
\end{prop}
\proof{
 We consider constraint Equation \eqref{Bconstraint}.
Given the operators
$$ B(x,y) = \sum_{i,j=0}^{n-1} x w^{2ik} e_{i,i} \tp e_{j-l,j} + y w^{2jk} e_{i+l,i} \tp e_{j,j}, $$ 
and
$$ \tilde{p}^{\alpha} = \sum_{i,j = 0}^{n-1}[ w^{2jb} e_{i+a+j,i+a} \tp e_{i+j,i}
+ w^{-2bj} e_{i-a+j,i-a} \tp e_{i+j,i}], $$
we first calculate
\begin{align*}
	\tilde{p}^{(a,b)}& B(x,y)\\
	& = \sum_{i,j,s,t=0}^{n-1} [ w^{-2bj} e_{i-a+j,i-a} \tp  e_{i+j,i} + w^{2bj} e_{i+a+j,i+a} \tp e_{i+j,i}]
		[ x w^{2sk}e_{s,s} \tp e_{t-l,t} + yw^{2tk} e_{s+l,s} \tp e_{t,t}] \\
	& = \sum_{j,t=0}^{n-1} 
		[ xw^{2((t-l-a)k-bj)} e_{t-l-a+j,t-l-a} \tp e_{t-l+j,t} 
		+ xw^{2((t-l+a)k+bj)} e_{t-l+a+j,t-l+a} \tp e_{t-l+j,t}\\
	&  \hspace{12mm} + yw^{2(tk-bj)} e_{t-a+j,t-a-l} \tp e_{t+j,t} 
		+ yw^{2(tk+bj)} e_{t+a+j,t+a-l} \tp e_{t+j,t}].
\end{align*}
This leads to
\begin{eqnarray*}
	&&\tilde{p}^{(a,b)} B(x,y) \tilde{p}^{(c,d)} \\
	& = & \sum_{j,t=0}^{n-1} 
		[ xw^{2((t-l-a)k-bj)} e_{t-l-a+j,t-l-a} \tp e_{t-l+j,t} +
		xw^{2((t-l+a)k+bj)} e_{t-l+a+j,t-l+a} \tp e_{t-l+j,t}\\
	& & \hspace{12mm} + yw^{2(tk-bj)} e_{t-a+j,t-a-l} \tp e_{t+j,t} 
		+ yw^{2(tk+bj)} e_{t+a+j,t+a-l} \tp e_{t+j,t}]\tilde{p}^{(c,d)}\\
	& = & \sum_{j,t,v=0}^{n-1} 
	[     xw^{2((t-l-a)k-bj-dv)} \dn_{a+l}^{c} e_{t-l-a+j,t-c-v} \tp e_{t-l+j,t-v} \\
	& & \hspace{12mm} + xw^{2((t-l+a)k+bj-dv)} \dn_{c}^{l-a} e_{t-l+a+j,t-c-v} \tp e_{t-l+j,t-v} \\
	& & \hspace{12mm} + yw^{2(tk-bj-dv)} \dn_{c}^{a+l} e_{t-a+j,t-c-v} \tp e_{t+j,t-v} 
		+ yw^{2(tk+bj-dv)} \dn_{c}^{l-a} e_{t+a+j,t-c-v} \tp e_{t+j,t-v} \\
	& & \hspace{12mm} + xw^{2((t-l-a)k-bj+dv)} \dn_{-c}^{l+a} e_{t-l-a+j,t+c-v} \tp e_{t-l+j,t-v} \\
	& &	\hspace{12mm} + xw^{2((t-l+a)k+bj+dv)} \dn_{c}^{a-l} e_{t-l+a+j,t+c-v} \tp e_{t-l+j,t-v}\\
	& & \hspace{12mm} + yw^{2(tk-bj+dv)} \dn_{-c}^{a+l} e_{t-a+j,t+c-v} \tp e_{t+j,t-v} 
		+ yw^{2(tk+bj+dv)} \dn_{c}^{a-l} e_{t+a+j,t+c-v} \tp e_{t+j,t-v}]\\
	& = & n\sum_{j,v=0}^{n-1} 
	\biggl\{     w^{2(dv-bj)} \left[xw^{2((-l-a)k-bl)} + y\right] \dn_{k+b}^{d} \dn_{a+l}^{c} e_{j-a,v-c} \tp e_{j,v} \\
	& & \hspace{12mm}+ w^{2(dv+bj)}\left[xw^{2((-l+a)k+bl)} + y\right] \dn_{k-b}^{d} \dn_{c}^{l-a} e_{j+a,v-c} \tp e_{j,v} \\
	& & \hspace{12mm}+ w^{-(2dv+bj)}\left[xw^{2((-l-a)k-bl)} + y\right] \dn_{k+b}^{-d} \dn_{-c}^{l+a} e_{j-a,v+c} \tp e_{j,v} \\
	& & \hspace{12mm}+ w^{-(2dv-bj)}\left[xw^{2((-l+a)k+bl)} + y\right] \dn_{k-b}^{-d} \dn_{c}^{a-l}
e_{j+a,v+c} \tp e_{j,v} \biggr\}.
\end{eqnarray*}
We now substitute the above result into Equation (\ref{Bconstraint}),
$$p^{\alpha}[f_{\alpha}(xy^{-1})B(x,y) - f_{\beta}(xy^{-1}) B(y,x)]p^{\beta} = 0,$$
and finally arrive at the desired four constraints.

Very similar calculations for Equation \eqref{Cconstraint} leads to the same four constraint equations.}

\begin{lem} \label{lem1}
Let $w^2$ be a primitive $n$th root of unity and $l, k$ be integers such that $gcd(l,n)=gcd(k,n)=1$. Then
$$ \prod_{j=1}^{a} \left( \frac{z + w^{2l((2j-1)k + b)}}{1 + zw^{2l((2j-1)k + b)}} \right) 
= \prod_{j=1}^{\overline{a}} \left( \frac{z + w^{2l((2j-1)k + b)}}{1 + zw^{2l((2j-1)k + b)}} \right), $$
for all $a\in \N$, $b\in \Z.$
\end{lem}
The proof is omitted.  Henceforth we choose $l,k$ satisfying $gcd(l,n)=gcd(k,n)=1$ so that we may use the above lemma.
\begin{prop} \label{S'}
Let $S' = \{ (a,b) | a \geq 0, b \in \Z \}$, and consider a set of functions $\{f_{(a,b)}| (a,b) \in S'\}$. 
If the functions satisfy the relations  
\begin{eqnarray*}
f_{(a+l,b+k)}(z) & = \es & \left( 
\frac{z + w^{2((a+l)k+bl)}}{1 + zw^{2((a+l)k+bl)}} \right) f_{(a,b)}(z), \\
f_{(0,b)}(z) &=  \es & f_{(0,-b)}(z), \\
f_{(a,b)}(z) & = \es & f_{(a,b+n)}(z), \\
f_{(a,b)}(z) & = \es & f_{(a+n,b)}(z),  \\
f_{(\overline{a},b)}(z) & = \es & f_{(n-\overline{a},-b)}(z), \  \ (a,b) \in S',
\end{eqnarray*}

\noindent then the functions also satisfy the constraints given in Proposition \ref{prop3-4} for all 
$(a,b) \in S$.  Moreover, every set of functions satisfying the conditions of Proposition \ref{prop3-4} can 
be extended in a unique way to a set of functions defined on $S'$ satisfying the above conditions.  Hence 
the two sets of constraints are equivalent.
\end{prop}
The proof is straightforward and omitted.
\begin{lem} \label{lem2}
The set of functions
$$ 
f_{(a,b)}(z) = \prod_{j=1}^{\overline{al^{-1}}} 
\left( \frac{z + w^{2l((2j-1)k + b-akl^{-1})}}{1 + zw^{2l((2j-1)k + b-akl^{-1})}} \right) 
f_{(0,b-akl^{-1})}(z) 
$$
satisfies the conditions in Proposition \ref{S'} given
$$ f_{(0,b)}(z) = f_{(0,-b)}(z) = f_{(0,b+n)}, $$
for $b \in \Z$.
\end{lem}
The proof follows from Lemma \ref{lem1}.
Substituting the functions given in Lemma \ref{lem2} and their associated projection operators into our form for $\check{R}(z)$ gives the operator
$$ \check{R}(z) = \frac{1}{n}\sum_{i,j,a,b=0}^{n-1} w^{2bj} \prod_{p=1}^{\overline{al^{-1}}} 
\left( \frac{z + w^{2l((2p-1)k+b-akl^{-1})}}{1 + zw^{2l((2p-1)k+b-akl^{-1})}} \right) 
f_{(0,b-akl^{-1})}(z) e_{i+a+j,i+a} \tp e_{i+j,i}, $$
which satisfies Equation $(\ref{eqnLLRc})$.
Note that in these calculations the only property that we have used is that $w^{2}$ 
is a primitive $n$th root of unity, irrespective of whether $n$ is odd or even. 
This allows us to use these calculations later on without alteration.

For $\check{R}(z)$ to be a descendent, it must also satisfy the Yang--Baxter equation. 
For the moment we will set $l=k=1$, and consider the more general case in Section \ref{other reps}. Then the functions become
\begin{equation} 
f_{(a,a + b)}(z) = \prod_{j=1}^{a} 
\left( \frac{z + w^{2(2j-1+b)}}{1 + zw^{2(2j-1+b)}} \right) f_{(0,b)}(z). 
\label{ph1}
\end{equation}
This gives the operator
\begin{eqnarray*}
\check{R}(z) 
&= & f_{(0,0)}(z)\left \{ \sum_{a=0}^{\frac{n-1}{2}} \prod_{j=1}^{a} 
\left( \frac{z + w^{2(2j-1)}}{1 + zw^{2(2j-1)}} \right) p^{(a,a)} \right \} \\
&  & + \sum_{b=1}^{\frac{n-1}{2}} f_{(0,b)}(z) \left\{p^{(0,b)} 
+ \sum_{a=1}^{\frac{n}{2}-1}\left[\prod_{j=1}^{a} 
\left( \frac{z + w^{2(2j-1 + b)}}{1 + zw^{2(2j-1+ b)}} \right) 
p^{(a,\overline{a+b})} \right. \right. \\
&  & \left. \left. + \prod_{j=1}^{a} 
\left( \frac{z + w^{2(2j-1 - b)}}{1 + zw^{2(2j-1-b)}} \right) 
p^{(a,\overline{a-b})} \right] \right\}, 
\end{eqnarray*}
or equivalently
\begin{equation}
\check{R}(z) = \sum_{i,j,a=0}^{n-1} 
\left[\frac{1}{n} \sum_{b=0}^{n-1} w^{2bj} \prod_{p=1}^{a} 
\left( \frac{z + w^{2(2p-1+b-a)}}{1 + zw^{2(2p-1+b-a)}} \right) 
f_{(0,b-a)}(z) \right] e_{i+a+j,i+a} \tp e_{i+j,i}. 
\label{ph2}
\end{equation}
In $\check{R}(z)$ there remain $\frac{n+1}{2}$ arbitrary functions, one for each disconnected component in 
the tensor product graph of Figure \ref{FigTensorOdd}. This is in some contrast to the examples studied in \cite{zgb91,dgz94,dgz96} for which the tensor product graphs are connected.
To ensure $R(z)$ inherits properties from the representation of the canonical element, 
we enforce additional conditions. Firstly, we impose that $R(z)$ is self-adjoint, 
which is equivalent to
$$ \check{R}_{12}(z) = \check{R}_{21}^{\dagger}(z), \hspace{1cm} z \in \R,$$
where $\dagger$ is the adjoint operator. 

\begin{prop} \label{propad}
The following three statements about the matrix $\check{R}(z)$ of Equation (\ref{ph2})
and coefficient functions of Equation (\ref{ph1}) are equivalent:
\begin{itemize}
\item[(i)]
$\check{R}_{12}(z) = \check{R}_{21}^{\dagger}(z).$
\item[(ii)]
$ (f_{(a,b)}(z))^{*} = f_{(a,-b)}(z), \quad \forall (a,b) \in S'.$
\item[(iii)]
$(f_{(0,b)}(z))^{*} = f_{(0,b)}(z) = f_{(0,b+2c)}(z), \quad \forall z \in \R$ and $b,c \in \Z$.
\end{itemize}
Here $^{*}$ denotes complex conjugation.
\end{prop}
\proof{
Assume $(i)$ holds. Looking at the general form of the projection operators we find that
$$ (\tilde{p}^{(a,b)}_{21})^{\dagger} = (\tilde{p}^{(a,n-b)}_{12}). $$
Recalling the form of $\check{R}(z)$ given in Equation
(\ref{Rcheck}) and using the linear independence of the projection operators,
this implies that
$$ (f_{(a,b)}(z))^{*} = f_{(a,-b)}(z), \quad \forall (a,b) \in S'. $$

Conversely, $(i)$ follows directly from $(ii)$, and hence statements $(i)$ and $(ii)$ are equivalent.  

Now assume $(ii)$ holds. Equation (\ref{ph1}) gives us
$$ f_{(1,1+b)}(z) = \left( \frac{z + w^{2(b+1)}}{1 + zw^{2(b+1)}} \right) f_{(0,b)}(z), $$
and
$$ f_{(1,-1-b)}(z) = \left( \frac{z + w^{-2(b+1)}}{1 + zw^{-2(b+1)}} \right) f_{(0,-b-2)}(z). $$
We see that
$$ (f_{(1,-1-b)}(z))^{*} = \left( \frac{z + w^{2(b+1)}}{1 + zw^{2(b+1)}} \right) 
(f_{(0,-b-2)}(z))^{*} = \left( \frac{z + w^{2(b+1)}}{1 + zw^{2(b+1)}} \right) 
f_{(0,b+2)}(z), $$
which implies that
$$ f_{(0,b)}(z) = f_{(0,b+2)}(z). $$
Combining this with the constraint of $f_{(0,b)}(z)$ we find
$$ (f_{(0,b)}(z))^{*} = f_{(0,b)}(z) = f_{(0,b+2c)}(z), $$
for $b,c \in \Z$.

Conversely, suppose $(iii)$ holds. 
It follows that
$$ 
f_{(a,b)}(z) = \prod_{j=1}^{a} \left( \frac{z + w^{2(2j-1+b-a)}}{1 +
zw^{2(2j-1+b-a)}} \right) f_{(0,b-a)}(z). 
$$ 
Hence 
\begin{eqnarray*}
  (f_{(a,b)}(z))^{*}
  & = & \left(\prod_{j=1}^{a} \left( \frac{z + w^{2(2j-1+b-a)}}{1 + zw^{2(2j-1+b-a)}} \right) f_{(0,b-a)}(z)\right)^{*} \\
  & = & \prod_{j=1}^{a} \left( \frac{z + w^{-2(2j-1+b-a)}}{1 + zw^{-2(2j-1+b-a)}} \right) \left(f_{(0,b-a)}(z)\right)^{*} \\
  & = & \prod_{j=1}^{a} \left( \frac{z + w^{2(2(-j)+1-b+a)}}{1 + zw^{2(2(-j)+1-b+a)}} \right) f_{(0,0-b-a)}(z) \\
  & = & \prod_{j=1}^{a} \left( \frac{z + w^{2(2j-1-b-a)}}{1 + zw^{2(2j-1-b-a)}} \right) f_{(0,-b-a)}(z) \\
  & = & f_{(a,-b)}(z). 
\end{eqnarray*} 

This completes the proof.
}

As the second index of the function can be considered modulo $n$ and $n$ is
odd, imposing self-adjointness implies that there is only one arbitrary function left. Without loss 
of generality we consider it to be $f_{(0,0)}(z)$. This can be seen as an 
overall scalar of our operator thus we are able to set it to a constant. Now we impose the limiting condition
\begin{equation} \label{limit}
\lim_{z \rightarrow 0}{R}(z) = \pm (\pi_n^{\pm} \tp \pi_n^{\pm}) \mathcal{R}.
\end{equation}

This sets the scalar to be $f_{(0,0)}(z) = 1$, so the functions become
$$ 
f_{(a,b)}(z)
 = \prod_{j=1}^{a} 
\left( \frac{z + w^{2(2j-1+b-a)}}{1 + zw^{2(2j-1+b-a)}} \right). 
$$
Hence 
$$ 
\check{R}(z)
 = \sum_{b=0}^{\frac{n-1}{2}} p^{(0,b)} + \sum_{a=1}^{\frac{n-1}{2}}\sum_{b=0}^{n-1} 
\prod_{j=1}^{a} \left( \frac{z + w^{2(2j-1+b-a)}}{1 + zw^{2(2j-1+b-a)}} \right) p^{(a,b)}, 
$$
or equivalently
$$ 
\check{R}(z) 
= \sum_{i,j,a=0}^{n-1} \left[\frac{1}{n} \sum_{b=0}^{n-1} w^{2bj} 
\prod_{p=1}^{a} \left( \frac{z + w^{2(2p-1+b-a)}}{1 + zw^{2(2p-1+b-a)}} \right) \right] 
e_{i+a+j,i+a} \tp e_{i+j,i}. 
$$
To simplify this we define the functions
$$ 
g_{(a,j)}(z) = \frac{1}{n} \sum_{b=0}^{n-1} w^{2bj} f_{(a,b)}(z)
 = \frac{1}{n} \sum_{b=0}^{n-1} w^{2bj} \prod_{p=1}^{a} 
\left( \frac{z + w^{2(2p-1+b-a)}}{1 + zw^{2(2p-1+b-a)}} \right) , 
$$
where $a \in \N$ and $j \in \Z$. In terms of these new functions $g_{(a,j)}(z)$, our 
operator becomes
\begin{equation*}
\check{R}(z) = \sum_{i,j,a=0}^{n-1} g_{(a,j)}(z) e_{i+a+j,i+a} \tp e_{i+j,i},
\end{equation*}

\noindent or

\begin{equation} \label{newR}
R(z) = \sum_{i,j,a=0}^{n-1} g_{(a,j)}(z) e_{i+j,i+a} \tp e_{i+a+j,i}.
\end{equation}

\begin{rmk}
The operator $R(z)$ of Equation (\ref{newR}) satisfies the following properties:

\begin{itemize}
\item[(i)] $R(z)^{*} = R(z), \; \forall z \in \mathbb{R}$,
\item[(ii)] $R^{t}(z) = R(z), \; \forall z \in \mathbb{C}$,
\item[(iii)] $R^{-1}(z) = R(z), \; \forall z \in \mathbb{C}$,
\item[(iv)] $R_{12}(z) R_{21}(z^{-1}) = I \tp I, \; \forall z \in \mathbb{C}$,
\item[(v)] $\underset{z \rightarrow 0}{\lim} R(z) = \pm(\pi_n^\pm \tp \pi_n^\pm) \mathcal{R}$,
\item[(vi)] $\underset{z \rightarrow 1}{\lim} R(z) = P$.

\end{itemize}
\end{rmk}

\begin{prop} \label{grel}
The operator $R(z)$ of Equation (\ref{newR}) is a descendant if and only if 
$$ 
\sum_{k=0}^{n-1} g_{(a,k-d)}(x)g_{(k,a-b)}(xy)g_{(b,c-k)}(y) 
= \sum_{k=0}^{n-1}  g_{(c,k-b)}(x)g_{(k,c-d)}(xy)g_{(d,a-k)}(y),
$$
for $0 \leq a,b,c,d \leq n-1$.
\end{prop}
This follows directly from the Yang--Baxter equation.

The matrix $R(z)$ of Equation (\ref{newR}) agrees with the $R$-matrix obtained in \cite{DIL2006} when $n=3$.  It has also been verified computationally for 
all odd $n \leq 17$ that the functions $g(z)$ satisfy the equation given in Proposition \ref{grel} above, and 
thus that $R(z)$ is a solution to the Yang--Baxter equation.  The computations were performed with 
$z$ treated as an arbitrary complex number and $w$ as an arbitrary primitive $n$th root of unity.

\begin{conj}
The matrix $R(z)$ given in Equation \eqref{newR} is a descendant of the zero-field six-vertex model with $D(D_n)$ symmetry for all odd $n$.
\end{conj}

\section{Descendants associated with $D(D_{2m})$}
In this section we construct a family of solutions of the Yang-Baxter equation
using $D(D_{2m})$, treating the cases where $m$ is even and
odd separately. It will be shown that although in both cases the 
representations have the same general form, the case where $m$ is odd has more in common 
with the solutions from the previous section associated with $D(D_n)$ where $n$ is
odd.

Throughout this entire section, we consider $w$ to be 
a primitive $2m$th root of unity and we use the $r(z)$ stated previously and 
associated with $\pi_{2} = \pi_{2}^{(k,l)}$.
The $L$-operator arising from the cases of $m$ odd or even can be considered at
once, since the approach is similar to that of Section \ref{ssecLopDnodd}.

\subsection{Construction of the $L$-operator}
In this subsection we construct an operator $L(z) \in {\rm End}(V_{2} \tp
V_{m})$ 
using a similar approach as before in Section \ref{ssecLopDnodd}. 
Again, we assume that the $L$-operator is of the form
$$ L(z) = (\pi_{2} \tp \pi_{m}) [R + h(z) (R^{T})^{-1}], $$
where $\pi_{m}$ is either $\pi_{m,\tau}^{(0,b)}$ or $\pi_{m,\sigma\tau}^{(0,b)}$ with $b\in \{0,1\}$. 
Applying this representation and an appropriate basis transformation on the two-dimensional space we find
$$
L(z) = \sum_{i=0}^{m-1}
\left \{ (w^{2ik}e_{1,2} + w^{-2ik}e_{2,1}) \tp e_{i,i} + h(z) 
\left[ e_{1,1} \tp e_{i-l,i} + e_{2,2} \tp e_{i+l,i} \right] \right\}. 
$$
We note that the basis transformation is of the same general form used 
in Section $\ref{ssecLopDnodd}$, as are the operators $r(z)$ and $L(z)$. Thus it follows that 
our $L$-operator is of the form
$$ 
L(z) 
= \sum_{i=0}^{m-1}
\left \{ (w^{2ik}e_{1,2} + w^{-2ik}e_{2,1}) \tp e_{i,i} + z 
\left[ e_{1,1} \tp e_{i-l,i} + e_{2,2} \tp e_{i+l,i} \right] \right\}. 
$$

\subsection{Construction of the descendants when $m$ is even}
\begin{figure}[ht]
\begin{center}
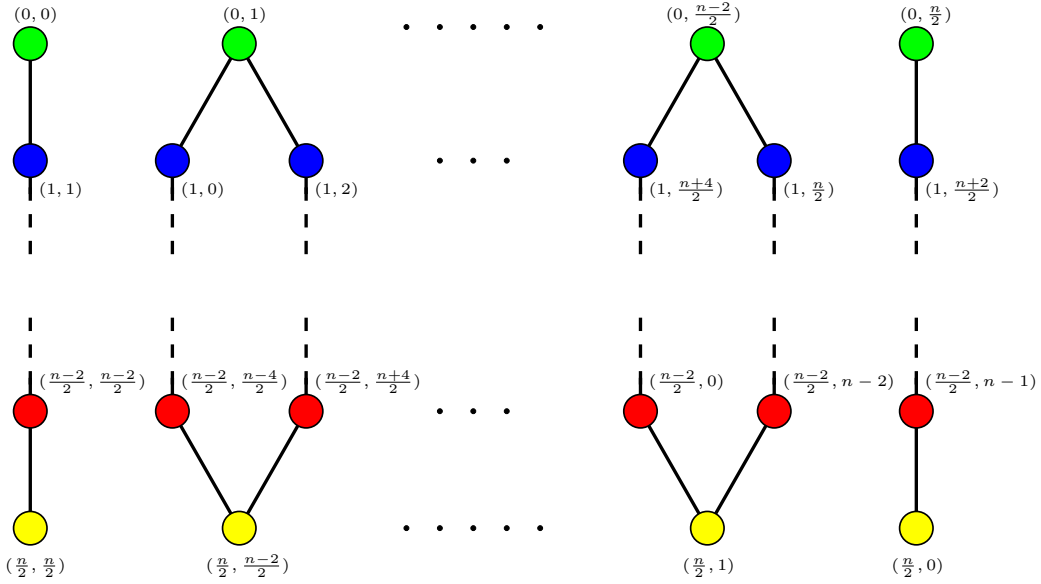
\end{center}
\caption{Tensor product graph for even $n$\label{FigTensorEven} when $l=k=1$.} 
\end{figure}

Here we consider the case of $D(D_{2m})$ where $m$ is even. We use $r(z)$
associated with $\pi_{2}^{(l,k)}$ with $w$ a primitive $2m$th root
of unity.
As before, we look for descendants of the form
$$ 
P R(z) = \check{R}(z) = \sum_{\alpha \in S} f_{\alpha}(z) p^{\alpha}, 
$$
where $p^{\alpha}$ are the projection operators previously calculated, 
$f_{\alpha}(z)$ are continuous functions and $S$ is the set of ordered pairs 
which correspond to non-zero projection operators. We require that $\check{R}(z)$ 
satisfies Equation $(\ref{eqnLLRc})$. As before, we use rescaled projection 
operators:
$$ 
\tilde{p}^{\alpha} 
= \sum_{i,j = 1}^{m}[ w^{2jb} e_{i+a+j,i+a} \tp e_{i+j,i} + w^{-2bj} e_{i-a+j,i-a} \tp e_{i+j,i}], \quad \alpha = (a,b) \in S.
$$
As these operators 
and the $L$-operator are equivalent to those found earlier, we can use the previous calculations to arrive at the tensor product graph shown in  
Figure $\ref{FigTensorEven}$.

We again define the set
$$ S' = \{(a,b) | a \in \N, b \in \Z \},$$
and extend the functions $f_\alpha(z)$ to that set in a way analogous to that of Proposition \ref{S'}.  We restrict ourselves 
to $k,l$ satisfying gcd$(l,m) = $gcd$(k,m) = 1$. Following the calculations of the previous section, we determine that if 
 $R(z)$ is a descendant then the functions $f_{(a,b)}(z)$ must satisfy
\begin{align*} 
f_{(a,akl^{-1} + b)}(z) 
&= \prod_{j=1}^{\overline{al^{-1}}} 
\left( \frac{z + w^{2l((2j-1)k + b)}}{1 + zw^{2l((2j-1)k + b)}} \right) f_{(0,b)}(z), \\
	f_{(0,b)}(z) & =  f_{(0,b+m)}(z), \\
	f_{(0,b)}(z) & =  f_{(0,-b)}(z), \qquad \forall (a,b) \in S'.
\end{align*} 

We now choose to set $l=k=1$; this gives the functions
$$ 
f_{(a,a + b)}(z) 
= \prod_{j=1}^{a} \left( \frac{z + w^{2(2j-1+b)}}{1 + zw^{2(2j-1+b)}} \right) f_{(0,b)}(z). 
$$
This produces the operator
\begin{eqnarray*}
\check{R}(z) 
&=& f_{(0,0)}(z)\left \{ \sum_{a=0}^{\frac{m}{2}} 
    \prod_{j=1}^{a} \left( \frac{z + w^{2(2j-1)}}{1 + zw^{2(2j-1)}} \right) p^{(a,a)} \right \} \\
& & + f_{(0,\frac{m}{2})}(z)\left \{ \sum_{a=0}^{\frac{m}{2}} \prod_{j=1}^{a} 
    \left( \frac{z + w^{2(2j-1 + \frac{m}{2})}}{1 + zw^{2(2j-1 + \frac{m}{2})}} \right) p^{(a,\overline{\frac{m}{2} + a})} \right\} \\
& & + \sum_{b=1}^{\frac{m}{2}-1} f_{(0,b)}(z) 
    \left\{p^{(0,b)} + \sum_{a=1}^{\frac{m}{2}-1}\left[\prod_{j=1}^{a} 
    \left( \frac{z + w^{2(2j-1+b)}}{1 + zw^{2(2j-1 + b)}} \right) p^{(a,\overline{a+b})} \right. \right. \\
& & \left. \left. + \prod_{j=1}^{a} 
    \left( \frac{z + w^{2(2j-1 - b)}}{1 + zw^{2(2j-1 - b)}} \right) 
    p^{(a,\overline{a-b})} \right] + \prod_{j=1}^{a} 
    \left( \frac{z + w^{2(2j-1+ \frac{m}{2})}}{1 + zw^{2(2j-1 + \frac{m}{2})}} \right) p^{(a,\overline{a+\frac{m}{2}})} \right\} \\
\end{eqnarray*}
or equivalently
$$ 
R(z) = \sum_{i,j,a=0}^{m-1} 
\left[\frac{1}{n}\sum_{b=0}^{m-1} w^{2bj} \prod_{p=1}^{a} 
\left( \frac{z + w^{2(2p-1+b-a)}}{1 + zw^{2(2p-1+ b-a)}} \right) 
f_{(0,b-a)}(z) \right] e_{i+j,i+a} \tp e_{i+a+j,i}. 
$$
Here we again enforce that our operator is self-adjoint. We are able to use the 
previous calculations and recall that
$$ 
(f_{(0,b)}(z))^{*} = f_{(0,b)}(z) = f_{(0,b+2c)}(z), 
$$
$\forall z \in \R$ and $b,c \in \N$. The second index of the function can be 
considered modulo $m$ and as $m$ is even the functions are partitioned into 
two sets. We see that we have only two functions in which we have any freedom 
left; without loss of generality we consider them to be $f_{(0,0)}(z)$ and 
$f_{(0,1)}(z)$. As there are two functions we find that we cannot consider them 
an overall scalar, so we need to impose additional conditions. As the $m$-dimensional 
representations of the canonical element are self-adjoint, we enforce that 
$R(z)$ is unitary . From this it follows that
$$ 
f_{(0,0)}(z) = \pm 1 \mbox{ and } f_{(0,1)}(z) = \pm 1. 
$$
Imposing the limiting condition given in Equation (\ref{limit}) sets
$$ 
f_{(0,0)}(z) = f_{(0,1)}(z) = 1. 
$$
This yields the operator
\begin{equation}
\check{R}(z) = \sum_{i,j,a=0}^{m-1} \left[\frac{1}{m} \sum_{b=0}^{m-1} w^{2bj} 
\prod_{p=1}^{a} \left( \frac{z + w^{2(2p-1+b-a)}}{1 + zw^{2(2p-1+b-a)}} \right) \right] 
e_{i+a+j,i+a} \tp e_{i+j,i}. 
\label{rcheckeven}
\end{equation}
To simplify this we define the functions
$$ g_{(a,j)}(z) = \frac{1}{m} \sum_{b=0}^{m-1} w^{2bj} f_{(a,b)}(z) 
= \frac{1}{m} \sum_{b=0}^{m-1} w^{2bj} \prod_{k=1}^{a} 
\left( \frac{1 + zw^{2(a-b+1-2k)}}{z + w^{2(a-b+1-2k)}} \right) , 
$$
where $a \in \N$ and $j \in \Z$. With the use of these new functions our operator 
becomes
\begin{equation} 
R(z) = \sum_{i,j,a=0}^{m-1} g_{(a,j)}(z) e_{i+j,i+a} \tp e_{i+a+j,i}. 
\label{Reven}
\end{equation}

\begin{rmk} The operator $R(z)$ of Equation (\ref{Reven}) satisfies the following properties:

\begin{itemize}
\item[(i)] $R(z)^{*} = R(z), \; \forall z \in \mathbb{R}$,
\item[(ii)] $R^{t}(z) = R(z), \; \forall z \in \mathbb{C}$,
\item[(iii)] $R^{-1}(z) = R(z), \; \forall z \in \mathbb{C}$,
\item[(iv)] $R_{12}(z) R_{21}(z^{-1}) = I \tp I, \; \forall z \in \mathbb{C}$,
\item[(iv)] $\underset{z \rightarrow 0}{\lim} R(z) = \pm(\pi_n^\pm \tp \pi_n^\pm) \mathcal{R}$,
\item[(v)] $\underset{z \rightarrow 1}{\lim} R(z) \neq P$.

\end{itemize}
\end{rmk}

Note that the last property shows that $R(z)$ does not satisfy regularity, unlike the $R(z)$ constructed from 
$D(D_n)$ when $n$ is odd.  
We do, however, have the following analogue of Proposition \ref{grel}:
\begin{prop} \label{geven}
The operator of Equation (\ref{Reven}) is a descendant if and only if the following constraint 
is satisfied:
$$ 
\sum_{k=0}^{m-1} g_{(a,k-d)}(x)g_{(k,a-b)}(xy)g_{(b,c-k)}(y) 
= \sum_{k=0}^{m-1}  g_{(c,k-b)}(x)g_{(k,c-d)}(xy)g_{(d,a-k)}(y), 
$$
for $0 \leq a,b,c,d \leq m-1$. 
\end{prop}
It has been computationally verified that the functions $g(z)$ satisfy the conditions in 
Proposition \ref{geven} above for all even $m \leq 16$, and hence that $R(z)$ is a solution 
to the Yang--Baxter equation.

\begin{conj}
The matrix $R(z)$ given by Equation \eqref{Reven} is a descendent of the zero-field six-vertex 
model with $D(D_{2m})$ symmetry for all even $m$.
\end{conj}

\subsection{Construction of the descendants when $m$ is odd}
For completeness, we include the descendants associated with $D(D_{2m})$ for odd
$m$. In our construction of the descendant of $r(z)$ we use a linear 
combination of projection operators found in Subsection \ref{sub2-3-2}. 
Using these projections, for a self-adjoint descendant of $r(z)$ which limits to 
a representation of the canonical element, we obtain the operator
$$ 
R(z) = \sum_{i,j,a=0}^{m-1} \left[\frac{1}{m}\sum_{b=0}^{m-1} w^{bj} 
\prod_{p=1}^{a} \left( \frac{z + w^{2(2p-1+b-a)}}{1 + zw^{2(2p-1+ b-a)}} 
\right) \right] e_{i+j,i+a} \tp e_{i+a+j,i}. 
$$
This is equivalent to the operator $R(z)$ found using $D(D_{m})$ where $m$ is odd.  Hence it 
obeys the properties stated in the remark in Section \ref{props}. Moreover, $R(z)$ is a descendant 
of the zero-field six-vertex model with $D(D_{2m})$ symmetry for odd $m \leq 17$ and conjectured to be a descendant for all odd $m \geq 3$.

\section{Generalised descendants}
In this section we investigate some of the choices we made in order to obtain 
descendants. Specifically, we investigate the effects of choosing a different initial 2-dimensional 
irrep and of imposing fewer conditions on the descendant.

\subsection{Dependency on choice of irreducible representations}\label{other reps}
In earlier sections we constructed descendants by starting from the $R$-matrix associated with 
the two-dimensional irrep $\pi_2^{(1,1)}$. Here we explore the consequences of starting with a different 
two-dimensional irrep, namely any $\pi_2^{(l,k)}$ where $\gcd(l,n) = \gcd(k,n)=1.$  We consider 
$n$-dimensional descendants for any $n \geq 3$, 
with $w$ a primitive $n$th root of unity. We recall that 
$$ r(z) = 
	\left(
	\begin{array}{cccc}
		w^{kl}z^{-1} - w^{-kl}z & 0 & 0 & 0\\
		0 & z^{-1} - z & w^{kl} - w^{-kl} & 0\\
		0 & w^{kl} - w^{-kl} & z^{-1} - z & 0\\
		0 & 0 & 0 & w^{kl}z^{-1} - w^{-kl}z
	\end{array}
\right)$$
and
$$
L(z) = \sum_{i=0}^{n-1}\left \{ (w^{ik}e_{1,2} + w^{-ik}e_{2,1}) \tp e_{i,i} 
+ z \left[ e_{1,1} \tp e_{i-l,i} + e_{2,2} \tp e_{i+l,i} \right] \right\}
$$
where $0 \leq l \leq \frac{n}{2}$ and  $0 \leq k \leq n-1$, although here $w$ differs slightly 
from that used in Section \ref{r}.  We now repeat our earlier calculations for 
our more general $l$ and $k$, obtaining 
$$ 
\check{R}(z) = \frac{1}{n}\sum_{i,j,a,b=0}^{n-1} w^{bjk} 
\prod_{p=1}^{\overline{al^{-1}}} 
\left( \frac{z + w^{lk(2p-1+b-al^{-1})}}{1 + zw^{lk(2p-1+b-al^{-1})}} \right) 
f_{(0,bk-akl^{-1})}(z) e_{i+a+j,i+a} \tp e_{i+j,i}. 
$$
We consider the change of basis on the $n$-dimensional space which yields
$$ e_{i,j} \rightarrow e_{il^{-1},jl^{-1}}, \hspace{1cm} i \in \Z. $$
Under this change of basis we find that our $L$-operator becomes
$$ 
L(z) = \sum_{i=0}^{n-1}\left \{ (w^{ilk}e_{1,2} + w^{-ilk}e_{2,1}) \tp e_{i,i} 
+ z \left[ e_{1,1} \tp e_{i-1,i} + e_{2,2} \tp e_{i+1,i} \right] \right\}, 
$$
while $\check{R}(z)$ becomes
$$ 
\check{R}(z) = \frac{1}{n}\sum_{i,j,a,b=0}^{n-1} w^{bjlk} 
\prod_{p=1}^{a} \left( \frac{z + w^{lk(2p-1+b-a)}}{1 + zw^{lk(2p-1+b-a)}} \right) 
f_{(0,(b-a)k)}(z) e_{i+a+j,i+a} \tp e_{i+j,i}. 
$$
From this we see that these different choices of initial two-dimensional representation
provide equivalent descendants. The different choices yield different basis 
transformations, a permutation on the arbitrary functions and a change of the 
root of unity, which must remain primitive. Thus any 2-dimensional 
irrep satisfying $\mbox{gcd}(l,n) = \mbox{gcd}(k,n)=1$ results in equivalent descendants. 
This is unsurprising as $\check{R}(z)$ is real whenever $z \in \mathbb{R}$, and choosing different 
$l,k$ effectively just changes the root of unity being used. 

\subsection{Descendants with an extra parameter associated with $D(D_{2m})$ when $m$ is even. \label{secOtherEvenDes}}
Here we return to $D(D_{2m})$ where $m$ is even and we find more descendants by 
imposing fewer constraints. We use the general form
$$ \check{R}(z) = \sum_{i,j,a=0}^{m-1} \left[\frac{1}{m} \sum_{b=0}^{m-1} w^{bj} 
\prod_{p=1}^{a} \left( \frac{z + w^{(2p-1+b-a)}}{1 + zw^{(2p-1+b-a)}} \right) 
f_{(0,b-a)}(z) \right] e_{i+a+j,i+a} \tp e_{i+j,i} $$
where $w$ is now a primitive $m$th root of unity.  We recall that previously we imposed that $R(z)$ is self-adjoint, unitary and obeys a limiting 
condition. If we ignore the limiting condition but enforce the other two conditions 
we find that 
$$ f_{(0,2b)}(z) = f_{(0,0)}(z) = \pm 1 \hspace{0.7cm} \mbox{and} \hspace{0.7cm} 
f_{(0,2b+1)}(z) = f_{(0,1)}(z) = \pm 1, $$
for all $b \in \Z$. Without loss of generality we can set 
$$ f_{(0,0)}(z) = 1. $$
This gives us two possibilities; the first choice is $f_{(0,1)}(z) = 1$, 
which yields
\begin{equation}
\check{R}^{+}(z) = \sum_{i,j,a=0}^{m-1} \left[\frac{1}{m} \sum_{b=0}^{m-1} w^{bj} 
\prod_{p=1}^{a} \left( \frac{z + w^{(2p-1+b-a)}}{1 + zw^{(2p-1+b-a)}} \right) \right] 
e_{i+a+j,i+a} \tp e_{i+j,i}. 
\label{r++}
\end{equation}
This operator corresponds to the conjectured descendant given in Equation \ref{rcheckeven}. The second 
option is $f_{(0,1)}(z) = -1$, which gives the operator
\begin{equation}
\check{R}^{-}(z) = \sum_{i,j,a=0}^{m-1} \left[\frac{1}{m} 
\sum_{b=0}^{m-1} (-1)^{a+b} w^{bj} \prod_{p=1}^{a} 
\left( \frac{z + w^{(2p-1+b-a)}}{1 + zw^{(2p-1+b-a)}} \right) \right] 
e_{i+a+j,i+a} \tp e_{i+j,i}. 
\label{r--}
\end{equation}
The matrices $\check{R}^+(z)$ and $\check{R}^{-}(z)$ are different; nonetheless they share many properties, 
including that they square to the identity and obey the unitarity condition. Moreover, we have the following 
proposition:
\begin{prop} \label{propr+r-}
The descendant $\check{R}^{-}(z)$ of Equation (\ref{r--}) satisfies the Yang--Baxter 
equation if and only if $\check{R}^{+}(z)$ given by Equation (\ref{r++}) does. 
\end{prop}
\proof{
We use the identity 
\begin{equation*}
\prod_{p=c+1}^{c+\frac{m}{2}} \left( \frac{z + w^{(2p-1+b-a)}}{1 + zw^{(2p-1+b-a)}} \right)
=(-1)^{a+b+\frac{m}{2}}.
\end{equation*}
This implies that
$$ 
\check{R}^{-}(z) = (-1)^{\frac{m}{2}}\frac{1}{m} \sum_{i,j,a,b=0}^{m-1} 
\left[ w^{bj} \prod_{p=1}^{a+\frac{m}{2}} 
\left( \frac{z + w^{(2p-1+b-a)}}{1 + zw^{(2p-1+b-a)}} \right)
\right] e_{i+a+j,i+a} \tp e_{i+j,i}.
$$
We now consider a basis transformation. Given any $\lambda \in \C$, 
there exists a basis transformation which yields
$$ 
e_{i,j} \rightarrow  \lambda^{\overline{i} - \overline{j}} e_{i,j}. 
$$
We choose $\lambda$ such that our operator becomes (after scaling)
$$ 
\check{R}^{-}(z) = \frac{1}{m} \sum_{i,j,a,b=0}^{m-1} \left[ w^{(b+\frac{m}{2})j} 
\prod_{p=1}^{a+\frac{m}{2}} \left( \frac{z + w^{(2p-1+b-a)}}{1 + zw^{(2p-1+b-a)}} \right) 
\right]  e_{i+a+j,i+a} \tp e_{i+j,i}. 
$$
We define the functions
$$
g_{(a,j)}(z) = \sum_{b=0}^{m-1} w^{bj} \prod_{p=1}^{a} 
\left( \frac{z + w^{(2p-1+b-a)}}{1 + zw^{(2p-1+b-a)}} \right),
$$
which allow to write the operators 
$$ 
\check{R}^{+}(z) = \sum_{i,j,a=0}^{m-1} g_{(a,j)}(z) e_{i+a+j,i+a} \tp e_{i+j,i} 
$$
and
$$ 
\check{R}^{-}(z) = \sum_{i,j,a=0}^{m-1} g_{(a+\frac{m}{2},j)}(z) e_{i+a+j,i+a} \tp e_{i+j,i}. 
$$  
Hence $\check{R}^{-}(z)$ differs from $\check{R}^{+}(z)$ only by a 
basis transformation and a permutation of the entries. We calculated previously 
that $\check{R}^{+}(z)$ satisfies the Yang--Baxter equation if and only if
$$ 
\sum_{k=0}^{m-1} \left[g_{(a,k-d)}(x) g_{(k,a-b)}(xy) g_{(b,c-k)}(y) 
-  g_{(c,k-b)}(x)g_{(k,c-d)}(xy)g_{(d,a-k)}(y)\right] =0, 
$$
$0 \leq a,b,c,d \leq m-1$. Similarly we have that $\check{R}^{-}(z)$ satisfies 
the Yang--Baxter equation if and only if
\fontsize{10}{12} $$ \sum_{k=0}^{m-1} 
\left[g_{(a+\frac{m}{2},k-d)}(x) g_{(k+\frac{m}{2},a-b)}(xy) g_{(b+\frac{m}{2},c-k)}(y) 
-  g_{(c+\frac{m}{2},k-b)}(x)g_{(k+\frac{m}{2},c-d)}(xy)g_{(d+\frac{m}{2},a-k)}(y)\right] 
= 0 , 
$$ \fontsize{12}{14}
$0 \leq a,b,c,d \leq m-1$. As we can consider each of the indices of the functions 
modulo $m$ we find that the two above constraints are equivalent, 
hence the result.
}
Thus we have found another family of conjectured descendants. We explain the 
existence of $\check{R}^{-}(z)$ by considering the values of functions associated with the irreps. 
To obtain $\check{R}^{+}(z)$ we simply set every function 
associated with irreps from the conjugacy class $\{ e \}$ equal to one, i.e. $f_{(0,b)}=1$ for $b \in \mathbb{Z}$. Conversely,  
$\check{R}^{-}(z)$ is obtained by setting every function associated with irreps from the 
conjugacy class $\{ \sigma^{m} \}$ equal to one, i.e. $f_{(\frac{m}{2},b)}(1)=1$ for $b \in \mathbb{Z}$. We recall that $e$ and $\sigma^{m}$ 
are the central elements in $D_{2m}$ when $m$ is even. It is in part due to this enlarged centre that we obtain more general solutions in 
this case.

It is also possible for us to ignore the unitary condition and only impose that $R(z)$ is self-adjoint.  
This leads to the constraints
$$ 
f_{(0,2b)}(z) = f_{(0,0)}(z) = (f_{(0,0)}(z))^{*} \hspace{0.7cm} \mbox{and} 
\hspace{0.7cm} f_{(0,2b+1)}(z) = f_{(0,1)}(z) = (f_{(0,1)}(z))^{*}, 
$$
for all $b \in \Z$ and $z \in \R$. Using these constraints we set
$$ f_{(0,0)}(z) = 1 + f(z) \hspace{0.7cm} \mbox{and} \hspace{0.7cm} 
f_{(0,1)}(z) = 1 - f(z), $$
where $f(z)$ is an arbitrary real function. This yields the operator
\begin{equation} 
\check{R}(z) = \check{R}^{+}(z) + f(z) \check{R}^{-}(z). 
\label{Rcheckz}
\end{equation}
This is invertible provided $$ f(z) \neq \pm 1. $$
The function $f(z)$ is equivalent to a second
parameter. That is, the operator
\begin{equation}
\check{R}(z,\mu) = \check{R}^{+}(z) + \mu \check{R}^{-}(z),
\label{Rcheckzmu}
\end{equation}
is invertible for $\mu\neq \pm 1$. Furthermore, we have the following:
\begin{prop}
The operator $\check{R}(z,\mu)$ given by Equation (\ref{Rcheckzmu})
satisfies 
\begin{equation}
\check{R}_{12}(x,\lambda)\check{R}_{23}(xy, \mu)\check{R}_{12}(y, \nu) 
= \check{R}_{23}(y, \nu)\check{R}_{12}(xy, \mu)\check{R}_{23}(x,\lambda), \label{eqnYBE2par}
\end{equation}
if and only if $\check{R}(z)$ given by Equation (\ref{Rcheckz}) is a solution to 
the Yang--Baxter equation with no constraints on $f(z)$. 
\end{prop}
\proof{
Consider a solution of the Yang-Baxter equation $R(z) \in \mbox{End}(V_{n} \tp V_{n})$ 
which contains an arbitrary function $f(z)$. Furthermore suppose every entry of 
$R(z)$ can be written as a polynomial in terms of $z$ and $f(z)$. Let
$$ \Omega = R_{12}(x)R_{13}(xy)R_{23}(y) - R_{23}(y)R_{13}(xy)R_{12}(x). $$
We observe that every entry of $\Omega$ is expressible as
$$ \sum_{i,j,k = 0}^{\infty} h_{ijk}^{l}(x,y) f^{i}(x)f^{j}(xy)f^{k}(y), $$
where $h_{ijk}^{l}(x,y)$ are polynomials in $x$ and $y$ 
and $l$ indexes the entry of $\Omega$. As $R(z)$ is a solution of the YBE, i.e. $\Omega = 0$, and 
$f(z)$ is an arbitrary function, we deduce that.
$$ h_{ijk}^{l}(x,y) = 0, \hspace{1cm} \forall x,y \in \C / \{ 0,1, \infty \} 
\hspace{0.1cm} \mbox{s.t.} \hspace{0.1cm} x \neq y. $$
Let $y_0 \in \C / \{ 0,1, \infty \}$. There are at most we have four values of $x$ for which 
$h_{ijk}^{l}(x,y_0)$ can be non-zero, but $h_{ijk}^{l}(x,y_0)$ is continuous in $x$. 
This means that
$$ 
h_{ijk}^{l}(x,y_0) = 0 \hspace{1cm} \forall x \in \C. 
$$
By symmetry
$$ 
h_{ijk}^{l}(x,y) = 0, \hspace{1cm} \forall (x,y) 
\in (\C \times \C) / (\{ 0,1, \infty \} \times \{ 0,1, \infty \}). 
$$
Thus there are at most 9 points in which $h_{ijk}^{l}(x,y)$ can be non-zero; however 
$h_{ijk}^{l}(x,y)$ is continuous in $x$ and $y$. Hence
$$ 
h_{ijk}^{l}(x,y) = 0, \hspace{1cm} \forall x,y \in \C. 
$$
Let $R(z,\mu)$ be the operator derived from $R(z)$ in which we have 
replaced the free function $f(z)$ with $\mu$. Every entry of $R(z)$ must 
be expressible as a polynomial in terms of $z$ and $\mu$. If we let
$$ \Omega = R_{12}(x,\lambda)R_{13}(xy, \mu)R_{23}(y, \nu) 
- R_{23}(y, \nu)R_{13}(xy, \mu)R_{12}(x,\lambda), $$
$x,y,\lambda,\mu,\nu \in \C$, then we find every entry of $\Omega$ can be 
written
$$ 
\sum_{i,j,k = 0}^{\infty} h_{ijk}^{l}(x,y) \lambda^{i}\mu^{j}\nu^{k}. 
$$
As shown previously
$$ h_{ijk}^{l}(x,y) = 0, \hspace{1cm} \forall x,y \in \C. $$
Thus $\Omega = 0$ and $R(z,\mu)$ satisfies
$$ R_{12}(x,\lambda)R_{13}(xy, \mu)R_{23}(y, \nu) 
= R_{23}(y, \nu)R_{13, \mu}(xy)R_{12}(x,\lambda). $$

Conversely, suppose we have $R(z,\mu)$ which satisfies
$$ 
R_{12}(x,\lambda)R_{13}(xy, \mu)R_{23}(y, \nu) 
= R_{23}(y, \nu)R_{13}(xy, \mu)R_{12}(x,\lambda). 
$$
If we consider $R(z) = R(z,f(z))$ where $f(z)$ is an arbitrary function then 
$R(z)$ must satisfy the YBE.
To recover the result, we use the fact that $\check{R}(z)=P\ R(z)$ and
$\check{R}(z,\mu)=P\ R(z,\mu)$. 
This proof works if the entries of $R(z)$ can be written as a 
polynomial $f(z)$ whose coefficients are rational functions of $z$. We are able 
to scale $R(z)$ by the product of all the denominators of the coefficients of 
$f(z)$, hence turning it into a polynomial.
}
We are able to verify using Maple that $\check{R}(z)$ given by Equation (\ref{Rcheckz}) 
satisfies the Yang--Baxter equation for even $n$ up to 12, and hence that 
$\check{R}(z,\mu)$ satisfies Equation (\ref{eqnYBE2par}). We conjecture that this holds 
true for all even $n$, and have shown it to hold in the limit $x=y=0$.

\section{Connection to the Fateev--Zamolodchikov model}
Closely associated with the Yang-Baxter equation is the star-triangle relation (STR), given by
$$ \sum_{d=0}^{N-1} \bar{W}(x|a-d) W(xy|d,c) \bar{W}(y|d-b) = W(x|b-c) \bar{W}(xy|a-b) W(y|a-c),$$
for $0 \leq a,b,c \leq N-1$. One well-known solution of the STR is the $N$-state Fateev--Zamolodchikov model \cite{FateevZam1982b}, which has weights 
$$ W(z|l) = \prod_{j=1}^{l} \frac{\lambda^{2j-1}z-1}{\lambda^{2j-1}-z} \hspace{0.7cm} \mbox{and} \hspace{0.7cm} \bar{W}(z|l) = \prod_{j=1}^{l} \frac{\lambda^{2j-1} - \lambda z}{\lambda^{2j}z  - 1}, $$
for $0 \leq l \leq N-1$, where $\lambda$ is a primitive $2N$th root of unity. These weights are extended by the relations
$$ W(z|l) = W(z|N+l) \hspace{0.7cm} \mbox{and} \hspace{0.7cm} \bar{W}(z|l) = \bar{W}(z|N+l),  $$
while also satisfying
$$ W(z|l) = W(z|N-l) \hspace{0.7cm} \mbox{and} \hspace{0.7cm} \bar{W}(z|l) = \bar{W}(z|N-l), $$
for $l \in \Z$. These weights satisfy the STR and lead to the $R$-matrix defined by
$$ R(\tilde{x},\tilde{y}) = \sum_{a_{1},a_{2},b_{1},b_{2}=1}^{N} R_{a_{1}a_{2}}^{b_{1}b_{2}}(\tilde{x},\tilde{y}) e_{b_{1},a_{1}} \tp e_{b_{2},a_{2}}, $$
where
$$ \tilde{x} ={x_{1}\choose x_{2}}, \hspace{1cm} \tilde{y} = {y_{1}\choose y_{2}} $$
and
$$ R_{a_{1}a_{2}}^{b_{1}b_{2}}(\tilde{x},\tilde{y}) = \bar{W}(x_{1}y_{1}^{-1}|a_{1}-b_{2}) W(x_{2}y_{1}^{-1}|a_{1}-a_{2}) \bar{W}(x_{2}y_{2}^{-1}|a_{2}-b_{1}) W(x_{1}y_{2}^{-1}|b_{2}-b_{1}). $$ 

Through a private communication with V. Bazhanov and J. Perk \cite{BazhPerk2009} we learnt of a connection between the $D(D_{n})$ solution and the Fateev--Zamolodchikov model. Specifically in the private communication a limiting case of the 3-state Fateev--Zamolodchikov model was shown to reduce to the $D(D_{3})$ solution. Using the ideas presented in \cite{BazhPerk2009} we are able to establish a connection between the $N$-state Fateev--Zamolodchikov model and $D(D_{n})$ (or equivalently $D(D_{2n})$) solution in the case where $N=n$ and $n$ is odd.

To investigate which limit of the Fateev--Zamolodchikov model leads to the $D(D_{n})$ model, we determine when the Fateev--Zamolodchikov $R$-matrix squares to the identity. Using certain properties of the weights it is possible to show that
$$ R(\tilde{x},\tilde{y})R(\tilde{x}^{-T},\tilde{y}^{-T}) \propto I \tp I, $$
where 
$$ \tilde{x} = {x_{1} \choose x_{2}}
\hspace{0.5cm} \mbox{and}\hspace{0.5cm} \tilde{x}^{-T} = {x_{2}^{-1} \choose x_{1}^{-1}}. $$
Thus the inverse of the $R$-matrix is known up to a scalar multiple. Furthermore we are able to show that if the $R$-matrix squares to the identity then it is equivalent to one which satisfies the  constraint
$$ \tilde{x} = \tilde{x}^{-T} \hspace{0.5cm} \mbox{and}\hspace{0.5cm} \tilde{y} = \tilde{y}^{-T}. $$

This implies that the $R$-matrix can be reduced to
$$ R(x,y) = \sum_{a_{1},a_{2},b_{1},b_{2}=1}^{N} R_{a_{1}a_{2}}^{b_{1}b_{2}}(x,y) e_{b_{1},a_{1}} \tp e_{b_{2},a_{2}}, $$
where
$$ R_{a_{1}a_{2}}^{b_{1}b_{2}}(x,y) = \bar{W}(xy^{-1}|a_{1}-b_{2}) W(x^{-1}y^{-1}|a_{1}-a_{2}) \bar{W}(x^{-1}y|a_{2}-b_{1}) W(xy|b_{2}-b_{1}). $$
To obtain the \textit{difference property} we set 
\begin{eqnarray}
 R(z) = \lim_{x,y\rightarrow \infty} R(x,y), 
 \label{fzz}
 \end{eqnarray}
where $z= \frac{x}{y}$. For odd $n \leq 11$ we are able to computationally verify that this $R$-matrix is equivalent (up to a basis transformation) to the $D(D_{n})$ $R$-matrix while setting $\lambda = - w^{-1}$. For odd $n>11$ we can verify that $R(0)$ is indeed equivalent to $(\pi_{n}^{+} \tp \pi_{n}^{+})\mathcal{R}$.

We now briefly comment on the descendants obtained from $D(D_{2n}),$ $n$ even.  There are two distinct $R$-matrices, $R^{+}(z)$ and $R^{-}(z)$, which both square to the identity.  
The multiplicities of the eigenvalues of $R^{+}(z)$ and $R^{-}(z)$ differ; furthermore neither eigenvalue spectrum matches that of $R(z)$ as defined by Equation (\ref{fzz}). 
Hence unlike the case when $n$ is odd, the $R$-matrix \eqref{Rcheckz} is not equivalent up to basis transformation of a limit of the Fateev--Zamolodchikov $R$-matrix.

\section{Summary}
In this paper we used the framework of descendants to construct $R$-matrices from the Drinfeld 
doubles of dihedral groups. 
For $3 \leq n \leq 17 $ and $w$ a primitive $n$th 
root of unity,
$$ 
R(z) = \sum_{i,j,a=0}^{n-1} \left[\frac{1}{n}\sum_{b=0}^{n-1} w^{bj} 
\prod_{p=1}^{a} \left( \frac{z + w^{(2p-1+b-a)}}{1 + zw^{(2p-1+ b-a)}} 
\right) \right] e_{i+j,i+a} \tp e_{i+a+j,i}. 
$$
is a descendant of the six-vertex model
$$
	r(z) = 
	\left(
	\begin{array}{cccc}
		wz^{-1} - w^{-1}z & 0 & 0 & 0\\
		0 & z^{-1} - z & w - w^{-1} & 0\\
		0 & w - w^{-1} & z^{-1} - z & 0\\
		0 & 0 & 0 & wz^{-1} - w^{-1}z
	\end{array}
	\right),
$$
with corresponding $L$-operator
$$ 
L(z) = \sum_{i=0}^{n-1}\left \{ (w^{i}e_{1,2} + w^{-i}e_{2,1}) \tp e_{i,i} 
+ z \left[ e_{1,1} \tp e_{i-1,i} + e_{2,2} \tp e_{i+1,i} \right] \right\}.
$$ 
We conjecture that this holds true for all $n \geq 2$.

We also showed that when $n$ is even we obtain $R$-matrices with a second, non-spectral parameter. 
Specifically, let $n>2$ be an even integer and $w$ a primitive $n$th root of unity. Given
$$ R^{+}(z) = \sum_{i,j,a=0}^{n-1} \left[\frac{1}{n} \sum_{b=0}^{n-1} w^{bj} 
\prod_{p=1}^{a} \left( \frac{z + w^{(2p-1+b-a)}}{1 + zw^{(2p-1+b-a)}} \right) 
\right] e_{i+j,i+a} \tp e_{i+a+j,i} $$
and
$$ R^{-}(z) = \sum_{i,j,a=0}^{n-1} \left[\frac{1}{n} \sum_{b=0}^{n-1} 
(-1)^{b-a} w^{bj} \prod_{p=1}^{a} \left( \frac{z + w^{(2p-1+b-a)}}{1 + zw^{(2p-1+b-a)}} 
\right) \right] e_{i+j,i+a} \tp e_{i+a+j,i}, $$
then
$$ R(z, \mu) = R^{+}(z) + \mu R^{-}(z) $$
satisfies Equation $(\ref{eqnYBE2par})$ for $n \leq 12$ and is conjectured to satisfy it for larger $n$. Moreover, 
$R(z, \mu)$ and $L(z)$ together satisfy
$$ 
L_{12}(x)L_{13}(y)R_{23}(x^{-1}y, \mu) = R_{23}(x^{-1}y, \mu)L_{13}(y)L_{12}(x). 
$$

{\bf Acknowledgements} The authors thank H. Au-Yang, V. Bazhanov, V. Mangazeev, and J. Perk for their insightful comments.  
P.S.I. acknowledges the support of an Early Career Researcher Grant from The University of Queensland.  K.A.D. acknowledges 
the support of the Australian Research Council under Discovery Project DP1092513.

\end{document}